\documentclass[11pt]{amsart}
\usepackage{amssymb,amsmath,amsthm}
\usepackage{graphicx,psfrag,color,pst-grad}
\usepackage{epsfig,graphicx,amssymb}

\newcommand{\be}{\begin{equation}}
\newcommand{\de}{\end{equation}}
\newcommand{\Frac}[2] {\frac{\textstyle #1} {\textstyle #2}}
\numberwithin{equation}{section}

\def\cal#1{\fam2#1}

\newtheorem{theorem}{Theorem}[section]
\newtheorem{definition}{Definition}[section]

\theoremstyle{definition}

\title[Climate dynamics: Natural variability and
related uncertainties] {Climate dynamics and fluid mechanics: \\
Natural variability and related uncertainties}
\author[M. Ghil, M.D. Chekroun, and E. Simonnet]{}

\email{ghil@lmd.ens.fr {\rm (M. Ghil)}}
 \email{chekro@lmd.ens.fr {\rm (M.D. Chekroun)}}
\email{eric.simonnet@inln.cnrs.fr {\rm (E. Simonnet)}}



\thanks{The present manuscript has been published as:
\textsc{M. Ghil, M.D. Chekroun and E.Simonnet}, {\it Physica D:
Nonlinear phenomena}, Special Issue on The Euler Equations: 250
Years On, {\bf 237}(2008), 2111--2126.}

\begin{document}
\maketitle

\vspace{-0.4cm}
 \centerline{\scshape Michael Ghil}
 {\footnotesize
 \centerline{D\'epartement Terre-Atmosph\`ere-Oc\'ean, Laboratoire de M\'et\'eorologie Dynamique (CNRS and
 IPSL) and}
   \centerline{Environmental Research and Teaching Institute,}
 \centerline{\'Ecole Normale Sup\'erieure, 75231 Paris Cedex 05, France;}
   \centerline{Dept. of Atmospheric Sciences and IGPP, University of California, Los Angeles, CA 90095-1565,
USA}} \vspace{1ex}

\centerline{\scshape Micka$\ddot{\mbox{e}}$l D. Chekroun}
{\footnotesize \centerline{Environmental Research and Teaching
Institute,} \centerline{\'Ecole Normale Sup\'erieure, 75231 Paris
Cedex 05, France}}

\vspace{1ex}

\centerline{\scshape Eric Simonnet} {\footnotesize
\centerline{Institut Non Lin\'eaire de Nice (INLN)-UNSA, UMR 6618
CNRS,} \centerline{1361, route des Lucioles 06560 Valbonne - France}
}


\begin{abstract}
The purpose of this review-and-research paper is twofold: (i) to
review the role played in climate dynamics by fluid-dynamical
models; and (ii) to contribute to the understanding and reduction of
the uncertainties in future climate-change projections. To
illustrate the first point, we review recent theoretical advances in
studying the wind-driven circulation of the oceans. In doing so, we
concentrate on the large-scale, wind-driven flow of the mid-latitude
oceans, which is dominated by the presence of a larger, anticyclonic
and a smaller, cyclonic gyre.  The two gyres share the eastward
extension of western boundary currents, such as the Gulf Stream or
Kuroshio, and are induced by the shear in the winds that cross the
respective ocean basins.  The boundary currents and eastward jets
carry substantial amounts of heat and momentum, and thus contribute
in a crucial way to Earth's climate, and to changes therein.

Changes in this double-gyre circulation occur from year to year and
decade to decade. We study this low-frequency variability of the
wind-driven, double-gyre circulation in mid-latitude ocean basins,
via the bifurcation sequence that leads from steady states through
periodic solutions and on to the chaotic, irregular flows documented
in the observations. This sequence involves local, pitchfork and
Hopf bifurcations, as well as global, homoclinic ones.

The natural climate variability induced by the low-frequency
variability of the ocean circulation is but one of the causes of
uncertainties in climate projections. The range of these
uncertainties has barely decreased, or even increased, over the last
three decades. Another major cause of such uncertainties could
reside in the structural instability---in the classical, topological
sense---of the equations governing climate dynamics, including but
not restricted to those of atmospheric and ocean dynamics.

We propose a novel approach to understand, and possibly reduce,
these uncertainties, based on the concepts and methods of random
dynamical systems theory. The idea is to compare the climate
simulations of distinct general circulation models (GCMs) used in
climate projections, by applying stochastic-conjugacy methods and
thus perform a stochastic classification of GCM families. This
approach is particularly appropriate given recent interest in
stochastic parametrization of subgrid-scale processes in GCMs.

As a very first step in this direction, we study the behavior of the
Arnol'd family of circle maps in the presence of noise. The maps'
fine-grained resonant landscape is smoothed by the noise, thus
permitting their coarse-grained classification.
\end{abstract}

\section{Introduction}

Charney {\it et al.} \cite{Charney} were the first to attempt a
consensus estimate of the equilibrium sensitivity of climate to
changes in atmospheric CO$_2$ concentrations. The result was the now
famous range for an increase of $1.5$ K to $4.5$ K in global
near-surface air temperatures, given a doubling of CO$_2$
concentration.

As the relatively new science of climate dynamics evolved through
the $1980$s and $1990$s, it became quite clear --- from
observational data, both instrumental and paleoclimatic, as well as
model studies
--- that Earth's climate never was and is unlikely to ever be in
equilibrium. The three successive IPCC reports (1991 \cite{hough91},
1996, and 2001 \cite{hough01}) concentrated therefore, in addition
to estimates of equilibrium sensitivity, on estimates of climate
change over the 21st century, based on several scenarios of CO$_2$
increase over this time interval, and using up to $18$ general
circulation models (GCMs) in the fourth IPCC Assessment Report (AR4)
\cite{SPM}.

The GCM results of temperature increase over the coming $100$ years
have stubbornly resisted any narrowing of the range of estimates,
with results for {\it T}$_s$ in $2100$ as low as $1.4$ K or as high
as $5.8$ K, according to the Third Assessment Report. The hope in
the research leading up to the AR4 was that a set of suitably
defined ``better GCMs" would exhibit a narrower range of year-$2100$
estimates, but this does not seem to have been the case.

The difficulty in narrowing the range of estimates for either
equilibrium sensitivity of climate or for end-of-the-century
temperatures is clearly connected to the complexity of the climate
system, the multiplicity and nonlinearity of the processes and
feedbacks it contains, and the obstacles to a faithful
representation of these processes and feedbacks in GCMs. The
practice of the science and engineering of GCMs over several decades
has amply demonstrated that any addition or change in the model's
``parametrizations" --- {\it i.e.}, of the representation of
subgrid-scale processes in terms of the model's explicit,
large-scale variables --- may result in noticeable changes in the
model solutions' behavior.

\begin{figure}[htpb]
\centering
\includegraphics*[width=9cm]{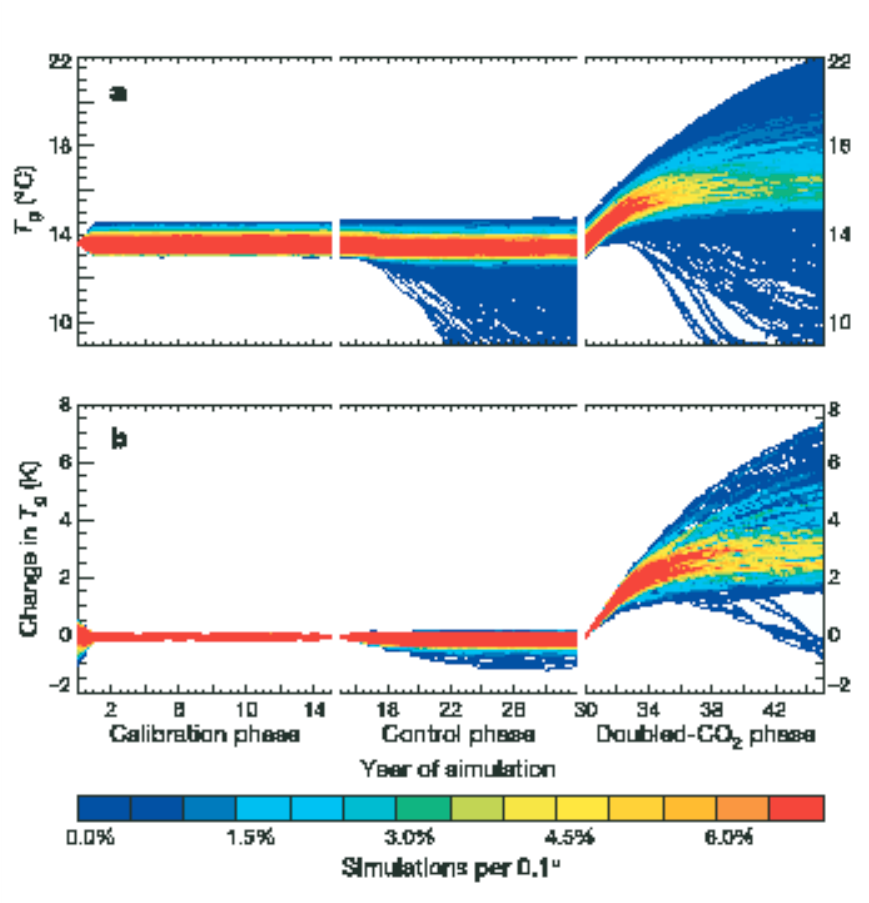}
\caption{Frequency distributions of global mean, annual mean,
near-surface temperature ($T_g$) for (a) 2,017 GCM simulations, and
doubled CO$_2$; and for (b) a subset of 414 stable simulations,
without substantial climate drift (from \cite{ltn}).}
\label{Stain_fig1}
\end{figure}

As an illustration, Fig. \ref{Stain_fig1} shows the sensitivity of
an atmospheric GCM, which does not include a dynamical ocean, to
changes in its model parameters. Several thousand simulations were
performed as part of the ``climate{\it prediction.}net" experiment
\cite{allen99}, using perturbations in several parameters of the
Hadley Centre's HadAM$3$ model \cite{pope}, coupled to a passive,
mixed-layer ocean model. The lower panel of Fig. \ref{Stain_fig1}
clearly illustrates a wide range of responses to CO$_2$ doubling,
from about $-1$ K to about $8$ K \cite{ltn}.

The last IPCC report \cite{SPM} has investigated climate change as a
result of various scenarios of CO$_2$ increase for a set of $18$
distinct GCMs. The best estimate of the temperature increase at the
end of the $21$st century from AR$4$ is about $4.0$$^\circ$ C for
the worst scenario of greenhouse-gas increase, namely A$1$F$1$, this
scenario envisages, roughly speaking, a future world with a very
rapid economic growth. The likely range of end-of-century increase
in global temperatures is of $2.4$--$6.4$$^\circ$ C in this case,
and comparably large ranges of uncertainties obtain for all the
other scenarios as well \cite{SPM}. The consequences of these
scentific uncertainties for the ethical quandaries arising in the
socio-economic and political decision-making process involved in
adaptation to and mitigation of climate changes are discussed in
\cite{RH_MG}.

An essential contributor to this range of uncertainty is natural
climate variability \cite{NRC95} of the coupled ocean-atmosphere
system. As mentioned already in \cite{MG0}, most GCM simulations do
not exhibit the observed interdecadal variability of the oceans'
buoyancy-driven, {\it thermohaline} circulation \cite{dijk_ghil}.
This circulation corresponds to a slow, pole-to-pole motion of the
oceans' main water masses, also referred to as the {\it overturning
circulation}. Cold and denser waters sink in the subpolar North
Atlantic and lighter waters rise over much wider areas of the lower
and southern latitudes.

Another striking example of low-frequency,
interannual-and-interdecadal variability is provided by the
near-surface, {\it wind-driven ocean circulation} \cite{dijk_ghil,
HDbook}. Key features of this circulation are described at length in
Section 2. The influence of strong thermal fronts
--- like the Gulf Stream in the North Atlantic or the Kuroshio in the North Pacific --- on the
mid-latitude atmosphere above is severely underestimated. Typical
spatial resolutions in the century-scale GCM simulations of
\cite{hough91, hough01, SPM, allen99,pope,ltn} are of the order of
$100$ km at best, whereas resolutions of $20$ km and less would be
needed to really capture the strong mid-latitude ocean-atmosphere
coupling just above the oceanic fronts \cite{FGS1,FGS2}.

An important additional source of uncertainty comes from the
difficulty to correctly parametrize global and regional effects of
clouds and their highly complex small-scale physics. This difficulty
is particularly critical in the tropics, where large-scale features
such as the El-Ni\~no/Southern Oscillation and the Madden-Julian
oscillation are strongly coupled with convective phenomena
\cite{MJO, Neel+98,GR00}.

The purpose of this paper is twofold. First, we describe in Section
2 the most recent theoretical results regarding the internal
variability of the mid-latitude wind-driven circulation, viewed as a
problem in nonlinear fluid mechanics. These results rely to a large
extent on the deterministic theory of dynamical systems
\cite{Arnold_geom, Guck}. Second, we address in Section 3 the more
general issue of uncertainties in climate change projections. Here
we rely on concepts and methods from random dynamical systems theory
\cite{LArnold} to help understand and possibly reduce these
uncertainties. Much of the material in the latter section is new; it
is supplemented by rigorous mathematical definitions and results in
Appendices A and B. A summary and an outlook on future work follow
in Section 4.


\section{Natural variability of the wind-driven ocean circulation}
\subsection{Observations}
To a first approximation, the main near-surface currents in the
oceans are driven by the mean effect of the winds. The trade winds
near the equator blow mainly from east to west and are called also
the tropical easterlies. In mid-latitudes, the dominant winds are
the prevailing westerlies, and towards the poles the winds are
easterly again. Three of the strongest near-surface,
mid-and-high-latitude currents are the Antarctic Circumpolar
Current, the Gulf Stream in the North Atlantic, and the Kuroshio
Extension off Japan. The Antarctic Circumpolar Current, sometimes
called the Westwind Drift, circles eastward around Antarctica; see
Fig. 2.
\begin{figure}[htpb]
\centering
\includegraphics*[width=8.5cm]{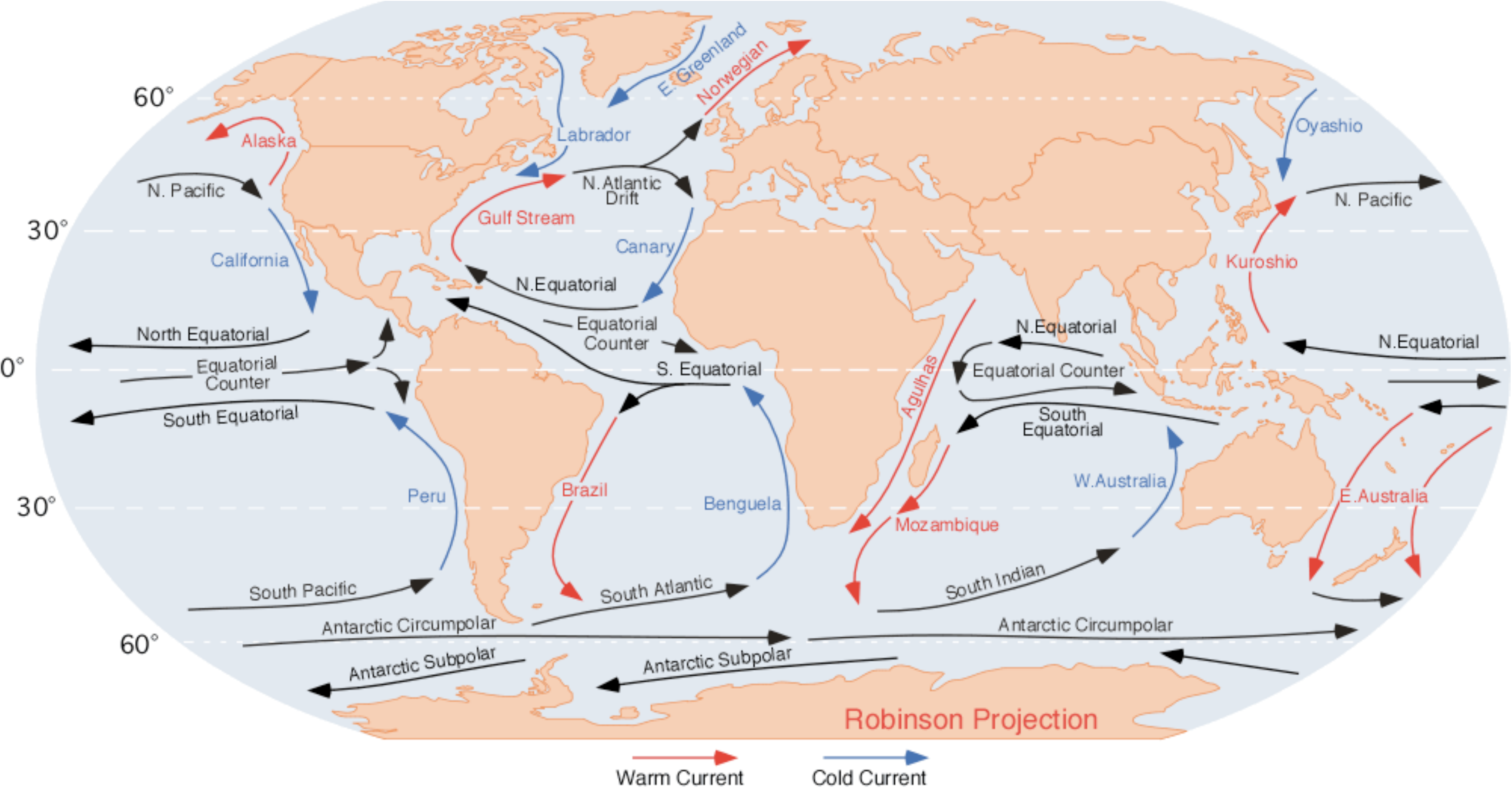}
\caption{A map of the main oceanic currents: warm currents in red
and cold ones in blue.} \label{oceanc}
\end{figure}

The Gulf Stream is an oceanic jet with a strong influence on the
climate of eastern North America and of western Europe. Actually,
the Gulf Stream is part of a larger, gyre-like current system, which
includes the North Atlantic Drift, the Canary Current and the North
Equatorial Current. It is also coupled with the pole-to-pole
overturning circulation. From Mexico's Yucatan Peninsula, the Gulf
Stream flows north through the Florida Straits and along the East
Coast of the United States. Near Cape Hatteras, it detaches from the
coast and begins to drift off into the North Atlantic towards the
Grand Banks near Newfoundland.

The Coriolis force is responsible for the so-called Ekman transport,
which deflects water masses orthogonally to the near-surface wind
direction and to the right \cite{Gill82,ghil-child,Ped87}. In the
North Atlantic, this Ekman transport creates a divergence and a
convergence of near-surface water masses, respectively, resulting in
the formation of two oceanic gyres: a smaller, cyclonic one in
subpolar latitudes, the other larger and anticyclonic in the
subtropics. This type of {\it double-gyre} circulation characterizes
all mid-latitude ocean basins, including the South Atlantic, as well
as the North and South Pacific.

The double-gyre circulation is intensified as the currents approach
the East Coast of North America due to the $\beta$-effect. This
effect arises primarily from the variation of the Coriolis force
with latitude, while the oceans' bottom topography also contributes
to it. The former, planetary $\beta$-effect is of crucial importance
in geophysical flows and induces free Rossby waves propagating
westward \cite{Gill82,ghil-child,Ped87}.

The currents along the western shores of the North Atlantic and of
the other mid-latitude ocean basins exhibit boundary-layer
characteristics and are commonly called western boundary currents
(WBCs). The northward-flowing Gulf Stream and the southward-flowing
Labrador Current extension meet near Cape Hatteras and yield  a
strong eastward jet. The formation of this jet and of the intense
recirculation vortices near the western boundary, to either side of
the jet, is mostly driven by internal, nonlinear effects.

\begin{figure}[htpb]
\centering
\includegraphics*[width=7.5cm]{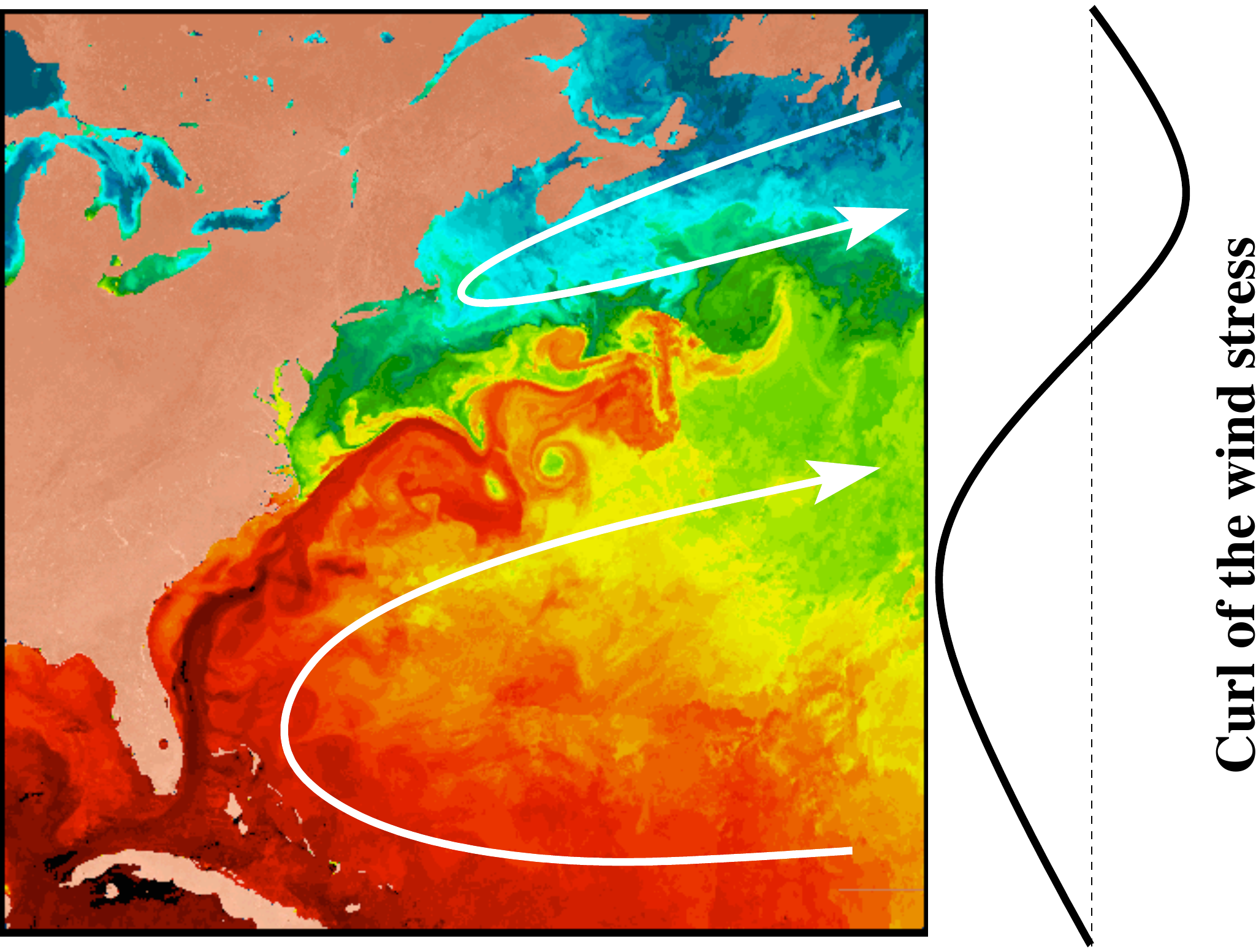}
 \caption{A satellite image of the sea surface temperature (SST)
over the northwestern North Atlantic (U.S. National Oceanic and
Atmospheric Administration), together with a sketch of the
associated double-gyre circulation. An idealized view of the amount
of potential vorticity injected into the ocean circulation by the
trade winds, westerlies and polar easterlies is shown to the right.}
\label{GS}
\end{figure}

Figure \ref{GS} illustrates how these large-scale wind-driven
oceanic flows self-organize, as well as the resulting eastward jet.
Different spatial and time scales contribute to this
self-organization, mesoscales eddies playing the role of the
synoptic-scale systems in the atmosphere. Warm and cold rings last
for  several months up to a year and have a size of about $100$ km;
two cold rings are clearly visible in Fig. \ref{GS}. Meanders
involve larger spatial scales, up to $1000$ km, and are associated
with interannual variability. The characteristic scale of the jet
and gyres is of several thousand kilometers and they exhibit their
own intrinsic dynamics on time scales of several years to possibly
several decades.

A striking feature of the wind-driven circulation is the existence
of two well-known North-Atlantic oscillations, with a period of
about $7$ and $14$ years, respectively. Data analysis of various
climatic variables, such as sea surface temperature (SST) over the
North Atlantic  or sea level pressure (SLP) over western Europe
\cite{Mor98,W99, DC&CDV} and local surface air temperatures in
Central England \cite{PGV}, as well as of proxy records, such as
tree rings in Britain, travertine concretions in southeastern France
\cite{Dubar}, and Nile floods over the last millennium or so
\cite{GK}, all exhibit strikingly robust oscillatory behavior with a
7-yr period and, to a lesser extent, with a 14-yr period. Variations
in the path and intensity of the Gulf Stream are most likely to
exert a major influence on the climate in this part of the world
\cite{HS65}. This is why theoretical studies of the low-frequency
variability of the double-gyre circulation are important.

Given the complexity of the processes involved, climate studies have
been most successful when using not just a single model but a full
hierarchy of models, from the simplest ``toy" models to the most
detailed GCMs \cite{GR00}. In the following, we describe one of the
simplest models of the hierarchy used in studying this problem.

\subsection{A simple model of the double-gyre circulation}
The simplest model that includes many of the mechanisms described
above is governed by the barotropic {\it quasi-geostrophic} (QG)
equations. The term geostrophic refers to the fact that large-scale
rotating flows tend to run parallel to, rather than perpendicular to
constant-pressure contours; in the oceans, these contours are
associated with the deviation from rest of the surfaces of equal
water mass, due to Ekman pumping. Geostrophic balance implies in
particular that the flow is divergence-free. The term barotropic, as
opposed to baroclinic, has a slightly different meaning in
geophysical fluid dynamics than in engineering fluid mechanics: it
means that the model describes a single fluid layer of constant
density and therefore the solutions do not depend on depth
\cite{Gill82, ghil-child, Ped87}.

We consider an idealized, rectangular basin geometry and simplified
forcing that mimics the distribution of vorticity contribution by
the winds, as sketched to the right of Fig. \ref{GS}. In our
idealized model, the amounts of subpolar and subtropical vorticity
injected into the basin are equal and the rectangular domain $\Omega
= (0,L_x) \times (0,L_y)$ is symmetric about the axis of zero wind
stress curl. The barotropic two-dimensional ($2$-D) QG equations in
this idealized setting are: \be \label{QGE}
\begin{array}{l}
q_t + J(\psi,q) -\nu \Delta^2 \psi + \mu \Delta \psi = -\tau \sin \Frac{2\pi y}{L_y} ,\\
q = \Delta \psi - \lambda_R^{-2} \psi + \beta y.
\end{array}
\de Here $q$ and $\psi$ are the potential vorticity and
streamfunction, respectively, and the Jacobian $J$ corresponds to
the advection of potential vorticity by the flow,  $J(\psi,q) =
\psi_x q_y - \psi_y q_x = {\bf u} \cdot \nabla q$, where ${\bf
u}=(-\psi_y,\psi_x)$, $x$ points east and $y$ points north. The
physical parameters are the strength of the planetary vorticity
gradient $\beta$, the Rossby radius of deformation $\lambda_R^{-2}$,
the eddy-viscosity coefficient $\nu$, the bottom friction
coefficient $\mu$, and the wind-stress intensity $\tau$. We use here
free-slip boundary conditions  $\psi = \Delta^2 \psi = 0$; the
qualitative results described below do not depend on the particular
choice of homogeneous boundary conditions.

We consider (\ref{QGE}) as an infinite-dimensional dynamical system
and study its bifurcation sets as the parameters change. Two key
parameters  are  the wind stress intensity $\tau$ and the eddy
viscosity $\nu$. An important property of (\ref{QGE}) is its mirror
symmetry in the $y=L_y/2$ axis. This symmetry can be expressed as
invariance with respect to the discrete $\mathbb{Z}_2$ group ${\cal
S}$: \be {\cal S}\left[ \psi(x,y) \right] = -\psi(x,L_y-y); \de any
solution of (\ref{QGE}) is thus accompanied by its mirror-conjugated
solution. Hence, in generic terms, the prevailing bifurcations are
of either the symmetry-breaking or the saddle-node or the Hopf type.

\subsection{Bifurcations in the  double-gyre problem}
The historical development of a comprehensive nonlinear theory of
the double-gyre circulation is interesting on its own, having seen
substantial progress in the last $15$ years. One can distinguish
four main steps.

\subsubsection{Symmetry-breaking bifurcations}
\begin{figure}[htpb]
\centering
\includegraphics*[width=9cm]{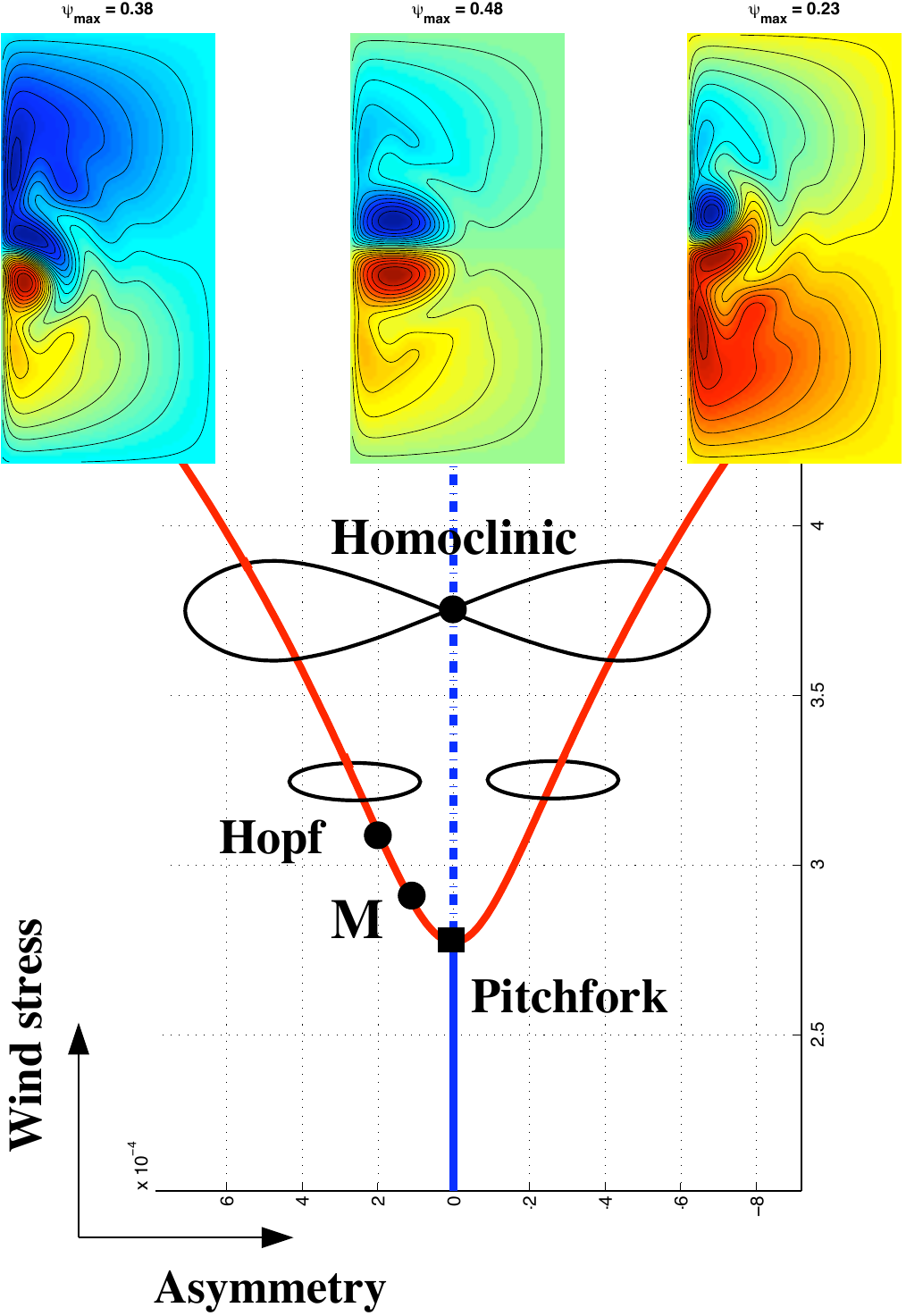}
\caption{Generic bifurcation diagram for the barotropic QG model of
the double-gyre problem: the asymmetry of the solution is plotted
versus the intensity of the wind stress $\tau$. The streamfunction
field is plotted for a steady-state solution associated with each of
the three branches; positive values in red and negative ones in blue
(after \cite{SGD05}).} \label{sk}
\end{figure}

The first step was to realize that the first generic bifurcation of
this QG model was a genuine pitchfork bifurcation that breaks the
system's symmetry as the nonlinearity becomes large enough
\cite{JJG0,JJG1,GC}. The situation is shown in Fig. \ref{sk}. When
the forcing is weak or the dissipation is large, there is only one
steady solution, which is antisymmetric with respect to the mid-axis
of the basin. This solution exhibits two large gyres, along with
their typical, $\beta$-induced WBCs.  Away from the western
boundary, such a near-linear solution (not shown) is dominated by
{\it Sverdrup} balance between wind stress curl and the meridional
mass transport \cite{Gill82,sv47}.

As the wind stress increases, the near-linear Sverdrup solution
develops an eastward jet  along the mid-axis, which penetrates
farther into the domain. This more intense, and hence more nonlinear
solution is still antisymmetric about the mid-axis, but loses its
stability for some critical value of the wind-stress intensity
(indicated by ``Pitchfork" in Fig. \ref{sk}).

A pair of mirror-symmetric solutions emerges and is characterized by
a rather different vorticity distribution; the streamfunction fields
associated with the two stable steady-state branches are plotted to
the upper-left and right of Fig. \ref{sk}. In particular, the jet in
such a solution exhibits a large meander, reminiscent of the one
seen in Fig. \ref{GS} just downstream of Cape Hatteras; note that
the colors in Fig. \ref{sk} have been chosen to facilitate the
comparison with Fig. \ref{GS}. These asymmetric flows are
characterized by one gyre being stronger in intensity than the other
and therefore the jet is deflected either to the southeast or to the
northeast.

\subsubsection{Gyre modes}
The next step was taken in part concurrently with \cite{JJG0, JJG1}
and in part shortly after \cite{SDG,DK,She} the first one. It
involved the study of time-periodic instabilities through Hopf
bifurcation from either an antisymmetric or an asymmetric steady
flow. Some of these studies concentrated on the wind-driven
circulation formulated for the stand-alone, single gyre \cite{She,
Ped96}. The idea was to develop a full generic picture of the
time-dependent behavior of the solutions in more turbulent regimes,
by classifying the various instabilities in a comprehensive way.
However, it quickly appeared that one kind of asymmetric
instabilities, called {\it gyre modes} \cite{JJG1,SDG}, was
prevalent across the full hierarchy of models of the double-gyre
circulation; furthermore, these instabilities trigger the lowest
nonzero frequency present in these models.

These modes always appear {\it after} the first pitchfork
bifurcation, and it took several years to really understand their
genesis: gyre modes arise as two eigenvalues merge --- one is
associated with a symmetric eigenfunction and responsible for the
pitchfork bifurcation, the other is associated with an antisymmetric
eigenfunction \cite{SD}; this merging is marked by $M$ in Fig.
\ref{sk}.

Such a phenomenon is not a bifurcation {\em stricto} {\em sensu}:
one has topological $C^0$ equivalence before and after the
eigenvalue merging, but not from the $C^1$ point of view. We recall
here that functions are $C^k$ if they and their inverses are $k$
times continuously differentiable. Still, this phenomenon is quite
common in small-dimensional dynamical systems with symmetry, as
exemplified by the unfolding of codimension-2 bifurcations of
Bogdanov-Takens type \cite{Guck}. In particular, the fact that gyre
modes trigger the lowest-frequency of the model is due to the
frequency of these modes growing quadratically from zero until
nonlinear saturation. Of course, these modes, in turn, become
unstable shortly after the merging, through a Hopf bifurcation,
indicated by ``Hopf" in Fig. \ref{sk}.

\subsubsection{Global bifurcations}
The importance of these gyre modes  was further confirmed recently
through an even more puzzling discovery. Several authors realized,
independently of each other, that the low-frequency dynamics of
their respective double-gyre models was driven by intense relaxation
oscillations of the jet \cite{SGITW0, Meach, Chang,NL, SGITW1,
SGITW2, SGD05}. These relaxation oscillations, already described in
\cite{JJG1, SDG}, were now attributed to {\it homoclinic}
bifurcations, with a global character in phase space
\cite{Guck,ghil-child}. In effect, the QG model reviewed here
undergoes a genuine homoclinic bifurcation (see Fig. \ref{sk}),
which is generic across the full hierarchy of double-gyre models.
Moreover, this global bifurcation is associated with chaotic
behavior of the flow due to the Shilnikov phenomenon
\cite{NL,SGD05}, which induces horseshoes in phase space.

The connection between such homoclinic bifurcations and gyre modes
was not immediately obvious, but Simonnet {\em et al.} \cite{SGD05}
emphasized that the two were part of a single, global dynamical
phenomenon. The homoclinic bifurcation indeed results from the
unfolding of the gyre modes' limit cycles. This familiar dynamical
scenario is again well illustrated by the unfolding of a
codimension-$2$ Bogdanov-Takens bifurcation, where the homoclinic
orbits emerge naturally. We deal, once more, with the
lowest-frequency modes, since homoclinic orbits have an infinite
period. Due to the genericity of this phenomenon, it was natural to
hypothesize that the gyre-mode mechanism, in this broader,
global-bifurcation context, gave rise to the observed $7$-yr and
$14$-yr North-Atlantic oscillations. Although this hypothesis may
appear a little farfetched, in view of the simplicity of the
double-gyre models analyzed in detail so far, it poses an
interesting question.

\subsubsection{Quantization and open questions}
The chaotic dynamics observed in the QG models after the homoclinic
bifurcation is eventually destroyed as the nonlinearity and the
resolution both increase. As one expects the real oceans to be in a
far more turbulent regime than those studied so far, some authors
proposed different mechanisms for low-frequency variability in fully
turbulent flow regimes \cite{Berloff, sk07}. It turns out, though,
that --- just as gyre modes could be reconciled with
homoclinic-driven dynamics,
--- the latter can also be reconciled with eddy-driven dynamics,
 via the so-called {\it quantization} of the low-frequency dynamics
\cite{S05}.

Primeau \cite{Primeau} showed that, in large basins comparable in
size with the North Atlantic, there is not only one but a set of
successive pitchfork bifurcations. One supercritical pitchfork
bifurcation, associated with the destabilization of antisymmetric
flows, is followed generically by a subcritical one, associated this
time with a stabilization of antisymmetric flows (modulo
high-frequency instabilities) \cite{S05}. As a matter of fact, this
phenomenon  appears to be a consequence of the spectral behavior of
the 2-D Euler  equations \cite{S07}, and hence of the closely
related barotropic  QG model in bounded domains.

Remarkably, this scenario repeats itself as the nonlinearity
increases, but now higher wavenumbers are involved in physical
space. Simonnet \cite{S05} showed that this was also the case for
gyre modes and the corresponding dynamics induced by global
bifurcations: the low-frequency dynamics is quantized as the jet
stream extends further eastward into the basin, due to the increased
forcing and nonlinearity.
\begin{figure}[htpb]
\centering
\includegraphics*[width=7cm]{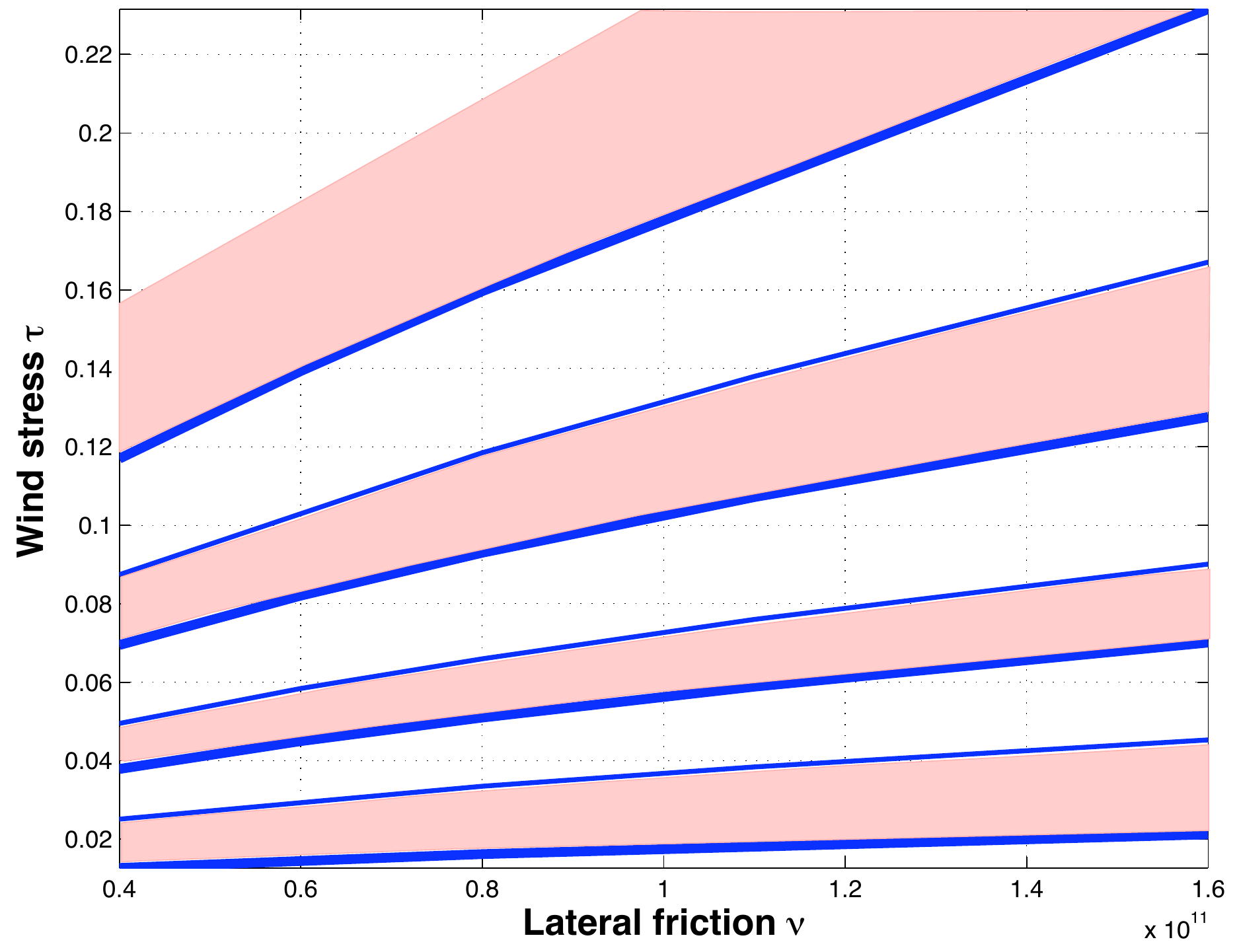}
\caption{Two-parameter plane, with the wind-stress intensity $\tau$
{\em vs.} the eddy-viscosity coefficient $\nu$: the curves indicate
the locations of supercritical and subcritical pitchfork
bifurcations. Each band is associated with a different wavenumber
and timescale (from \cite{S05}).} \label{Quanta}
\end{figure}
Figure \ref{Quanta} illustrates this situation: two families of
regimes can be identified, the colored bands correspond to
(supercritical) regimes driven by the gyre modes, the others to
(subcritical) regimes driven by the eddies. Note that this scenario
is also robust to perturbing the problem's symmetry.

The successive-bifurcation theory appears therewith to be fairly
complete for barotropic, single-layer models of the double-gyre
circulation. This theory also provides a self-consistent, plausible
explanation for the climatically important 7-year and 14-year
oscillations of the oceanic circulation and the related atmospheric
phenomena in and around the North-Atlantic basin
\cite{dijk_ghil,HDbook,Mor98,W99,DC&CDV,PGV,Dubar,GK,SGITW2,SGD05}.
The dominant 7- and 14-year modes of this theory also survive
perturbation by seasonal-cycle changes in the intensity and
meridional position of the westerly winds  \cite{Laxmi07}.

In baroclinic models, with two or more active layers of different
density, baroclinic instabilities \cite{dijk_ghil, FGS2,
Gill82,ghil-child, Ped87, HS65, Ped96, SGITW2, Berloff, sk07} surely
play a fundamental role, as they do in the observed dynamics of the
oceans. However, it is not known to what extent baroclinic
instabilities can destroy gyre-mode dynamics. The difficulty lies in
a deeper understanding of the so-called {\it rectification} process
\cite{DKrec}, which arises from the nonzero mean effect of the
baroclinic component of the flow.

Roughly speaking, rectification drives the dynamics far away from
any steady states. In this situation, dynamical systems theory
cannot be used as an explanation of complex, observed behavior
resulting from successive bifurcations that are rooted in a simple
steady state. Other tools from statistical mechanics and
nonequilibrium thermodynamics should, therefore, be considered
\cite{RSM,TT93, Bouchet1,Farrel}. Combining these tools with those
of the successive-bifurcation approach may eventually lead to a more
general and complete physical characterization of gyre modes in
realistic models.

\section{Climate-change projections and random dynamical systems (RDSs)}

As discussed in Section 1, the climate system's natural variability
and the difficulties in parametrizing subgrid-scale processes are
not the only causes for the uncertainties in projecting future
climate evolution. In this section, we address more generally these
uncertainties and present a novel approach for treating them. To do
so, we start with some simple ideas about deterministic {\em vs.}
stochastic modeling.

\subsection{Background and motivation}\label{3_1}

Many physical phenomena can be modeled by deterministic evolution
equations. Dynamical systems theory is essentially a geometric
approach for studying the asymptotic, long-term properties of
solutions to such equations in phase space. Pioneered by H.
Poincar\'e \cite{Poincare}, this theory took great strides over the
last fifty years. To apply the theory in a reliable manner to a set
of complex physical phenomena, one needs a criterion to evaluate the
{\it robustness} of a given model within a class of dynamical
systems. Such a criterion should help us deal with the inescapable
uncertainties in model formulation, whether due to incomplete
knowledge of the governing laws or inaccuracies in determining model
parameters.

In this context, Andronov and  Pontryagin \cite{AP} took a major
step toward classifying dynamical systems, by introducing the
concept of {\em structural stability}. Structural stability means
that a small, continuous perturbation of a given system preserves
its dynamics up to a {\it homeomorphism}, {\it i.e.}, up to a
one-to-one continuous change of variables that transforms the phase
portrait of our system into that of the nearby system; thus fixed
points go into fixed points, limit cycles into limit cycles, etc.
Closely related is the notion of {\it hyperbolicity} introduced by
Smale \cite{Smale}. A system is hyperbolic if, (very) loosely
speaking, its limit set can be continuously decomposed into
invariant sets that are either contracting or expanding; see
\cite{Katok} for more rigorous definitions.

A very simple example is the phase portrait in the neighborhood of a
fixed point of saddle type. In this case, the Hartman-Grobman
theorem states that the dynamics in this neighborhood is
structurally stable. The converse statement, {\it i.e.} whether
structural stability implies hyperbolicity, is still an open
question; the equivalence between structural stability and
hyperbolicity has only been shown in the $C^1$ case, under certain
technical conditions \cite{ro, rob, mane, Palis}. Bifurcation theory
is well grounded in the setting of hyperbolic dynamics. Problems
with hyperbolicity and bifurcations arise, however, when one deals
with more complicated limit sets.

Hyperbolicity was introduced initially to help pursue the
``dynamicist's dream'' of finding, in the abstract space of all
possible dynamical systems, an open and dense set consisting of
structurally stable ones. Being open and dense, roughly speaking,
means that any possible dynamical system can be approximated by
systems taken from this set, while systems in its complement are
negligible in a suitable sense.

Smale conjectured that hyperbolic systems form an open and dense set
in the space of all $C^1$ dynamical systems. If this conjecture were
true then hyperbolicity would be typical of all dynamics.
Unfortunately, though, this conjecture is only true for
one-dimensional dynamics and flows on disks and surfaces
\cite{peixoto}. Smale \cite{Smale66} himself found several
counterexamples to his conjecture. Newhouse \cite{Newhouse} was able
to generate open sets of nonhyperbolic diffeomorphisms using
homoclinic tangencies. For the physicist, it is even more striking
that the famous Lorenz attractor \cite{Lorenz63} is structurally
unstable. Families of Lorenz attractors, classified by topological
type, are not even countable \cite{GW, W}. In each of these
examples, we observe chaotic behavior in a nonhyperbolic situation,
{\em i.e., nonhyperbolic chaos.}

Nonhyperbolic chaos appears, therefore, to be a severe obstacle to
any ``easy" classification of dynamic behavior. As mentioned by
Palis \cite{Palis}, Kolmogorov already suggested at the end of the
sixties that ``the global study of dynamical systems could not go
very far without the use of new additional mathematical tools, like
probabilistic ones." Once more, Kolmogorov showed prophetic insight,
and nowadays the concept of {\it stochastic stability} is an
important tool in the study of genericity and robustness for
dynamical systems. To replace the failed program of classifying
dynamical systems based on structural stability and hyperbolicity,
Palis \cite{Palis} formulated the following {\it global conjecture}:
systems having only finitely many attractors ({\it i.e.} periodic or
chaotic sinks) -- such that (i) the union of their basins has full
Lebesgue measure; and (ii) each is stochastically stable in their
basins of attraction -- are dense in the $C^r, r \geq 1$ topology. A
system is stochastically stable if its Sinai-Ruelle-Bowen (SRB)
measure \cite{Sinai} is stable with respect to stochastic
perturbations, and the SRB measure is given by $\lim_{n\to \infty}
\frac{1}{n} \sum_i \delta_{z_i}$, with $z_i$ being the successive
iterates of the dynamics. This measure is obtained intuitively by
allowing the entire phase space to flow onto the attractor
\cite{ER85}.

Stochastic stability is fundamentally based on ergodic theory. We
would like to consider a more geometric approach, which can provide
a coarser, more robust classification of GCMs and their
climate-change projections. In this section, we propose such an
approach, based on concepts from the rapidly growing field of random
dynamical systems (RDSs), as developed by L. Arnold \cite{LArnold}
and his ``Bremen group," among others. RDS theory describes the
behavior of dynamical systems subject to external stochastic
forcing; its tools have been developed to help study the geometric
properties of stochastic differential equations (SDEs). In some
sense, RDS theory is the stochastic counterpart of the geometric
theory of ordinary differential equations (ODEs). This approach
provides a rigorous mathematical framework for a stochastic form of
robustness, while the more traditional, topological concepts do not
seem to be appropriate.

\subsection{RDSs, random attractors, and robust classification}
Stochastic parametrizations for GCMs aim at compensating for our
lack of detailed knowledge on small spatial scales in the best way
possible \cite{LN00,LN02,LN03,TP00,JPS,SZG}. The underlying
assumption is that the associated time scales are also much shorter
than the scales of interest and, therefore, the lag correlation of
the phenomena being parametrized is negligibly small. Stochastic
parametrizations thus essentially transform a deterministic
autonomous system into a nonautonomous one, subject to random
forcing.

Explicit time dependence in a dynamical system immediately raises a
technical difficulty. Indeed, the classical notion of attractor is
not always relevant, since any object in phase space is ``moving"
with time and the natural concept of forward asymptotics is
meaningless. One needs therefore another notion of attractor. In the
deterministic nonautonomous  framework, the appropriate notion is
that of a {\it pullback attractor} \cite{Langa}, which we present
below. The closely related notion of {\it random attractor} in the
stochastic framework is also explained briefly below, with further
details given in Appendix A.

\subsubsection{Framework and objectives}\label{Section_Amib_view_ofGC}

Before defining the notion of pullback attractor, let us recall some
basic facts about nonautonomous dynamical systems. Consider the ODE
\be \dot x = f(t,x) \label{ODE} \de on a vector space $X$; this
space could even be infinite-dimensional, if we were dealing with
partial or functional differential equations, as is often the case
in fluid-flow and climate problems. Rigorously speaking, we cannot
associate a dynamical system acting on $X$ with a nonautonomous ODE;
nevertheless, in the case of unique solvability of the initial-value
problem, we can introduce a two-parameter family of operators
$\{S(t,s)\}_{t \geq s}$ acting on $X$, with $s$ and $t$ real, such
that $S(t,s)x(s)=x(t)$ for $t\geq s$, where $x(t)$ is the solution
of the Cauchy problem with initial data $x(s)$. This family of
operators satisfies $S(s,s)=\mbox{Id}_X$ and $S(t,\tau) \circ
S(\tau, s)=S(t,s)$ for all $t\geq \tau \geq s$, and all real $s$.
This family of operators is called a ``process'' by Sell
\cite{sell}.It extends the classical notion of the resolvent of a
nonautonomous linear ODE to the nonlinear setting.

We can now define the pullback attractor as simply the family of
invariant sets $\{{\cal A}(t)\}$ that satisfy for every real $t$ and
all $x_0$ in $X$: \be \lim_{s \to -\infty} {\rm dist} ~(S(t,s) x_0,
{\cal A}(t)) = 0 \label{attr}. \de ``Pullback'' attraction does not
involve running time backwards; it corresponds instead to the idea
of measurements being performed at present time $t$ in an experiment
that was started at some time $s < t$ in the past: the experiment
has been running for long enough, and we are thus looking now at an
``attracting state." Note that there exists several ways of defining
a pullback attractor
--- the one retained here is a local one ({\it cf}. \cite{Langa} and
references therein);  see  \cite{berger} for further information on
nonautonomous dynamical systems in general.

In the stochastic context, noise forcing is modeled by a stationary
stochastic process. If the deterministic dynamical system of
interest is coupled to this stochastic process in a reasonable way
--- to be expressed below by the ``cocycle property"
--- then random pullback attractors may appear. These pullback
attractors will exist for almost each sample path of the driving
stochastic process, so that the same  probability distribution
governs both sample paths and their corresponding pullback
attractors. A more detailed explanation is given in Appendix A.

Roughly speaking, this concept of {\it random attractor} provides a
geometric framework for the description of asymptotic regimes in the
context of stochastic dynamics. To compare different stochastic
systems in terms of their random attractors that evolve in time, it
would be nice to be able to identify the common underlying geometric
structures via a random change of variables. This identification is
achieved through the concept of {\em stochastic equivalence} that is
developed in Appendix A, and it is central in obtaining a coarser
and more robust classification than in the purely deterministic
context.

Returning now to our main objective, suppose for instance that one
is presented with results from two distinct GCMs, say two
probability distributions functions (PDFs) of the temperature or
precipitation in a given area. These two PDFs are generated,
typically, by an ensemble of each GCM's simulations, as described in
the introduction, and they are likely to differ in their spatial
pattern. To ascertain the physical significance of this discrepancy,
one needs to know how each GCM result varies as either a
parametrization or a parameter value are changed.

In order to consider the difficult question of why GCM responses to
CO$_2$ doubling might differ, one idea is to investigate the
structure of the space of all GCMs. We mean therewith the space of
all deterministic GCMs, when their stochastic parametrizations are
switched off. We know, by now, from experience with GCM results over
several decades--- including the four IPCC assessment reports
\cite{hough91,hough01,SPM} and the climate{\it prediction}.net
exercise \cite{allen99,pope,ltn}---that there is enormous scatter in
this space; see also \cite{held,mcwill}. Our question, therefore,
is: can we achieve a more robust classification of GCMs when
stochastic parametrizations are used and for a given level of the
noise?

\begin{figure}[htpb]
\centering
\includegraphics*[width=5.5cm]{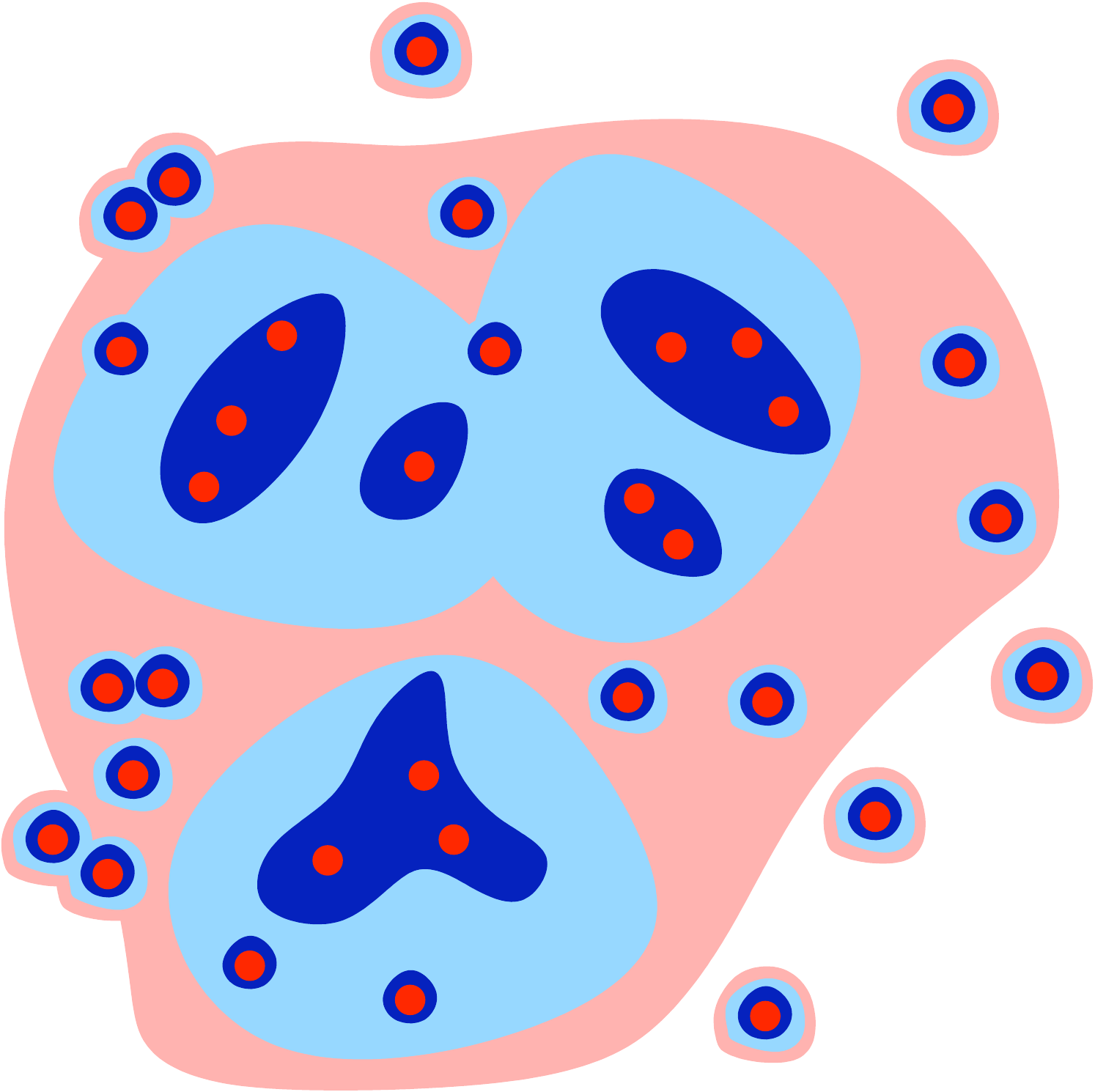}
\caption{A conjectural view of stochastic classification for GCMs,
using the concept of random attractors. Each point in red represents
a GCM in which stochastic parametrizations are switched off, while
each gray area represents a cluster of stochastically equivalent
GCMs for a given level of the noise.} \label{amib}
\end{figure}

As mentioned in Section \ref{3_1}, such a classification is not
feasible by restricting ourself to deterministic systems and
topological concepts. As one switches on stochastic parametrizations
\cite{LN00,LN02,LN03,TP00,JPS,SZG}, the situation might change, and
hopefully improve, dramatically: as the noise level becomes large
enough, the models' deterministic behavior may be completely
destroyed, and all the results could cluster into one huge, diffuse
clump. We would like, therefore, to investigate how a classification
based on stochastic equivalence evolves as the level of the noise or
the stochastic parametrizations change. As the noise tends to zero,
do we recover the ``granularity" of the set of all deterministic
dynamical systems? This idea is schematically represented in Fig.
\ref{amib}: for a given level of the noise, we expect the space of
all GCMs to be decomposed into a possibly finite number of classes.
Within one of these classes, all the GCMs are topologically
equivalent in the stochastic sense defined above; see Eq.
(\ref{stochequi}).

Serious difficulties might arise in this program, due to the
presence of nonhyperbolic chaos in climate models. Several studies
have pointed out that the characteristics of nonhyperbolic chaos in
the presence of noise may depend on its intensity and statistics
\cite{anischenko_et_al,anischenko_et_al2,
anischenko_et_al3,grebogi_et_al}.

Such issues, however, go well beyond the setting of this paper and
are left for further investigation. Much more modestly, we will
study here whether, in certain very simple cases, the conjectural
view of Fig. \ref{amib} might be relevant for some dynamical systems
that are ``metaphors" of climate dynamics. The following subsection
is dedicated to the study of such a metaphorical object, namely the
Arnol'd circle map.

\subsubsection{The stochastically perturbed circle map}
\label{Num_section}

To go beyond our pictorial view of stochastic classification for
GCMs in Fig. \ref{amib}, we study now the effect of noise on a
family of diffeomorphisms of the circle. This toy model exhibits two
features of interest for our purpose. The first one is that the
two-parameter family $\{F_{\tau,\epsilon}\}$ defined by Eq.
(\ref{amap}) below exhibits an infinite number of topological
classes \cite{Arnold_geom}. The second feature of interest is the
frequency-locking behavior observed in many field of physics in
general \cite{fl42,fl44, fl40} and in some
El-Ni\~no/Southern-Oscillation (ENSO) models in particular
\cite{EN1,EN2,EN3,EN4,EN6,SG01,EN5}. Studying noise effects on these
two features has, therefore, physical and mathematical, as well as
climatological relevance.

Many physical and biological systems exhibit interference effects
due to competing periodicities. One such effect is mode locking,
which is due to nonlinear interaction between an ``internal"
frequency $\omega_i$ of the system and an ``external" frequency
$\omega_e$. In the ENSO case, the external periodicity is the
seasonal cycle. A simple model for systems with two competing
periodicities is the well-known Arnol'd family of circle maps \be
\label{amap} x_{n+1} = F_{\tau,\epsilon}(x_n) := x_n + \tau -
\epsilon \sin(2\pi x_n)~{\rm mod}~1,  \de where basically
$\tau:=\omega_i/\omega_e$ and $\epsilon$ parameterizes the magnitude
of nonlinear effects; the map (\ref{amap}) is often called the {\it
standard circle map} \cite{Arnold_geom}.

These maps also represent frequency locking near a bifurcation of
Neimark-Sacker type  ({\it e.g.} \cite{Kuznetsov}, p. 434); here the
parameter $\tau$ is typically interpreted as the novel (internal)
frequency involved in the bifurcation and $\epsilon$ corresponds to
the nonlinearity near the bifurcation.

\begin{figure}[htpb]
\centering
\includegraphics*[width=10cm]{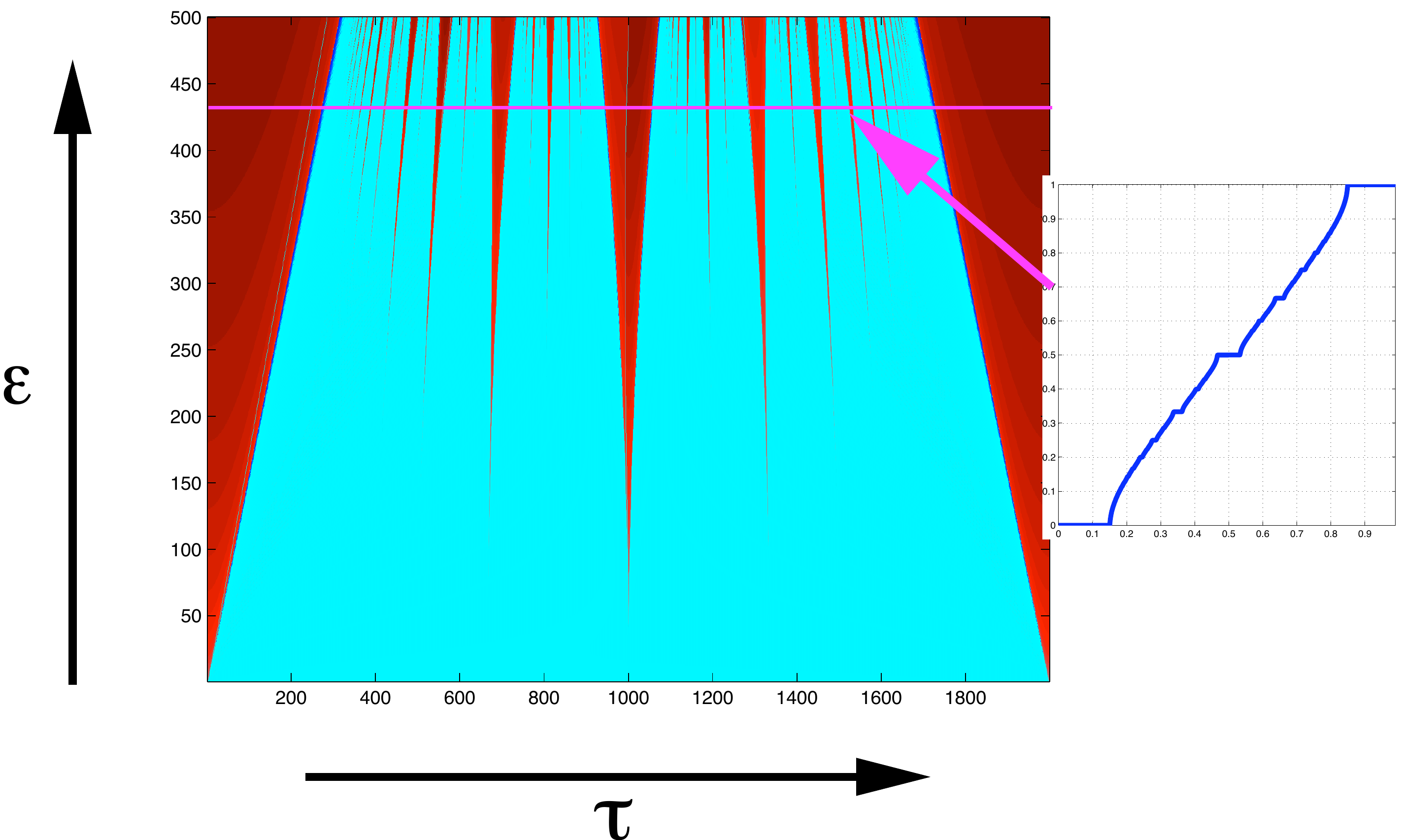}
\caption{Arnol'd tongues for the family of diffeomorphisms of the
circle; units for $\tau$ and $\epsilon$ are $5 \cdot 10^{-4}$ and
$10^{-4}$ respectively. Devil's staircase in the cross-section to
the right.} \label{tong}
\end{figure}

Such nonlinear coupling between two oscillators gives rise to a
characteristic pattern, in the plane of $\epsilon$ {\it vs.} $\tau$,
called Arnol'd tongues. We computed this pattern numerically for the
family of Eq. (\ref{amap}), together with a cross-section at a fixed
value of $\epsilon$; see Fig. \ref{tong}. This cross-section
exhibits the so-called Devil's staircase, with ``steps" on which the
{\it rotation number} \cite{Poincare} is constant within each
Arnol'd tongue; the rotation number measures the average rotation
per iterate of (\ref{amap}).

For $\epsilon=0$, two types of phenomena occur: either $\tau$ is
rational and in this case the dynamics is periodic with period $q$,
where $\tau = p/q$, or $\tau$ is irrational and the iterates
\{$x_n$\} fill the whole circle densely. As $\epsilon$ increases, an
Arnol'd tongue of increasing width grows out of each $\tau = p/q$ on
the abscissa $\epsilon = 0$. It follows that, in this very simple
case, such an Arnol'd tongue corresponds to hyperbolic dynamics that
is robust to perturbations, as verified by linearizing the map at
the periodic point; the rotation number is then rational and equal
to $p/q$.

The set of all these tongues is dense within the whole circle map
family, while the Lebesgue measure of this set, at given $\epsilon$,
tends to zero as $\epsilon$ goes to zero. On the contrary, if a
point in the $(\tau,\epsilon)$-plane does not belong to an Arnol'd
tongue, the rotation number for those parameter values is irrational
and the dynamics is nonhyperbolic; the latter fact follows, for
instance, from a theorem of Denjoy \cite{Denjoy} showing that such
dynamics is smoothly equivalent to an irrational rotation. The
probability to observe nonhyperbolic dynamics tends therewith to
unity as $\epsilon$ goes to zero. One has, therefore, a countably
infinite number of distinct topological classes, namely the Arnol'd
tongues $p/q$, and an uncoutably infinite number of maps with
irrational rotation numbers.

\begin{figure}[h]
\centering
\includegraphics*[width=10.5cm]{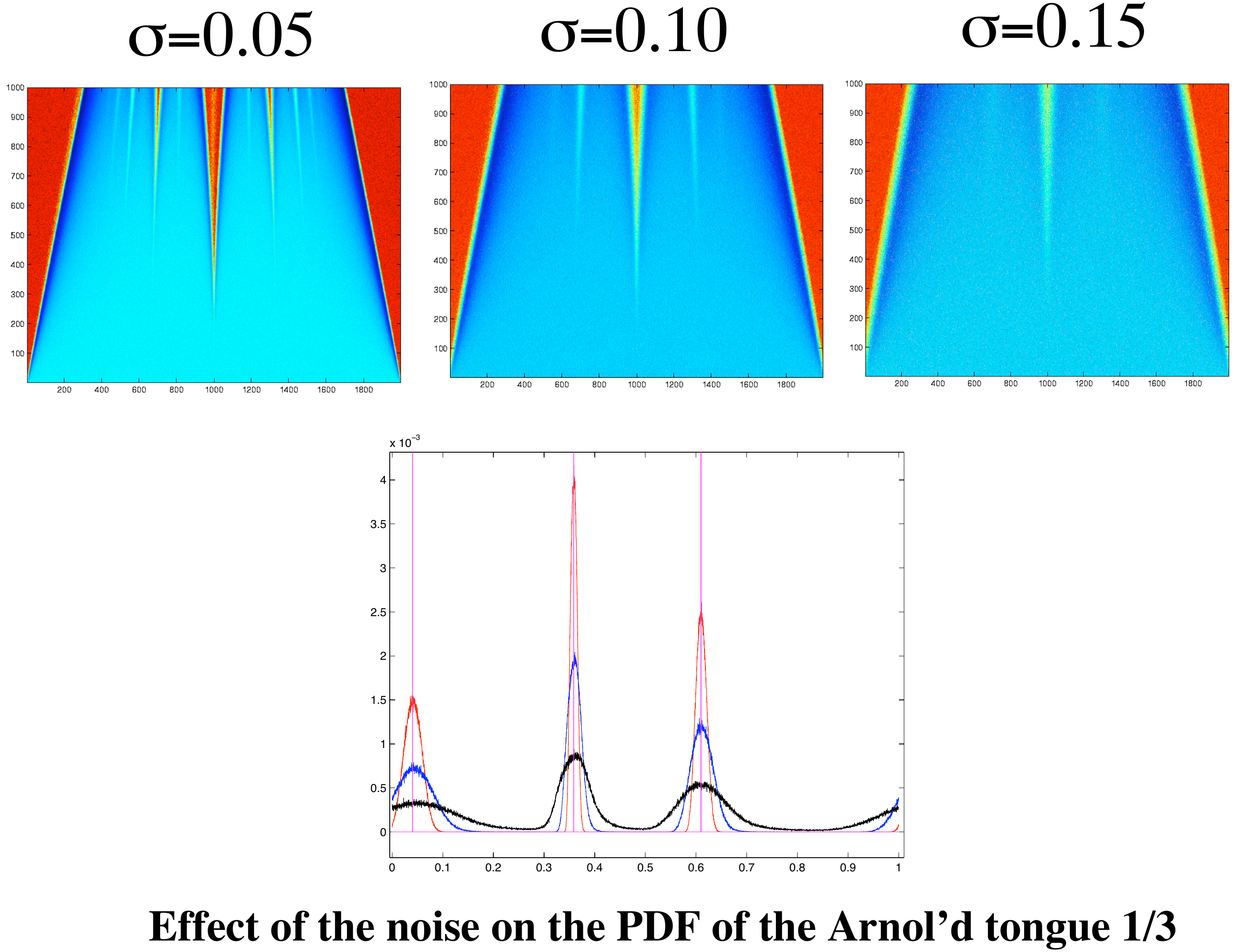}
\caption{Arnol'd tongues in the presence of additive noise with
different noise amplitudes $\sigma$. Upper panels: Arnol'd tongues
for $\sigma = 0.05, 0.10$ and $0.15$; lower panel: PDF for $\epsilon
= 0.9$ and the three $\sigma$-values in the upper panels: $\sigma =
0.05$ (red curve), $\sigma = 0.10$ (blue curve), and $\sigma = 0.15$
(black curve).} \label{noisytong}
\end{figure}

What happens when noise is added in Eq. (\ref{amap})? We consider
here the case of additive forcing by a noise process obtained via
sampling at each iterate $n$ a random variable with uniform density
and intensity $\sigma$. Experiments with colored, rather than white
noise and multiplicative, rather than additive noise led to the same
qualitative results. The results for additive white noise are shown
in Fig. \ref{noisytong} for three different levels of noise
intensity $\sigma$.

As expected, only the largest tongues survive the presence of the
noise; in particular, there is only a finite number of surviving
tongues, shown in red in Fig. \ref{noisytong}. Within such a
surviving tongue, the random attractor ${\cal A}(\omega)$ is a
random periodic cycle of period $q$ (not shown). In the blue region
outside the Arnol'd tongues, the random attractor is a fixed but
random point ${\cal A}(\omega) = \{ a(\omega) \}$: if one starts a
numerical simulation for a fixed realization of the noise $\omega$,
all initial data $x$ converge to the same fixed point $a$, say.

\begin{figure}[h]
\centering
\includegraphics*[width=9cm]{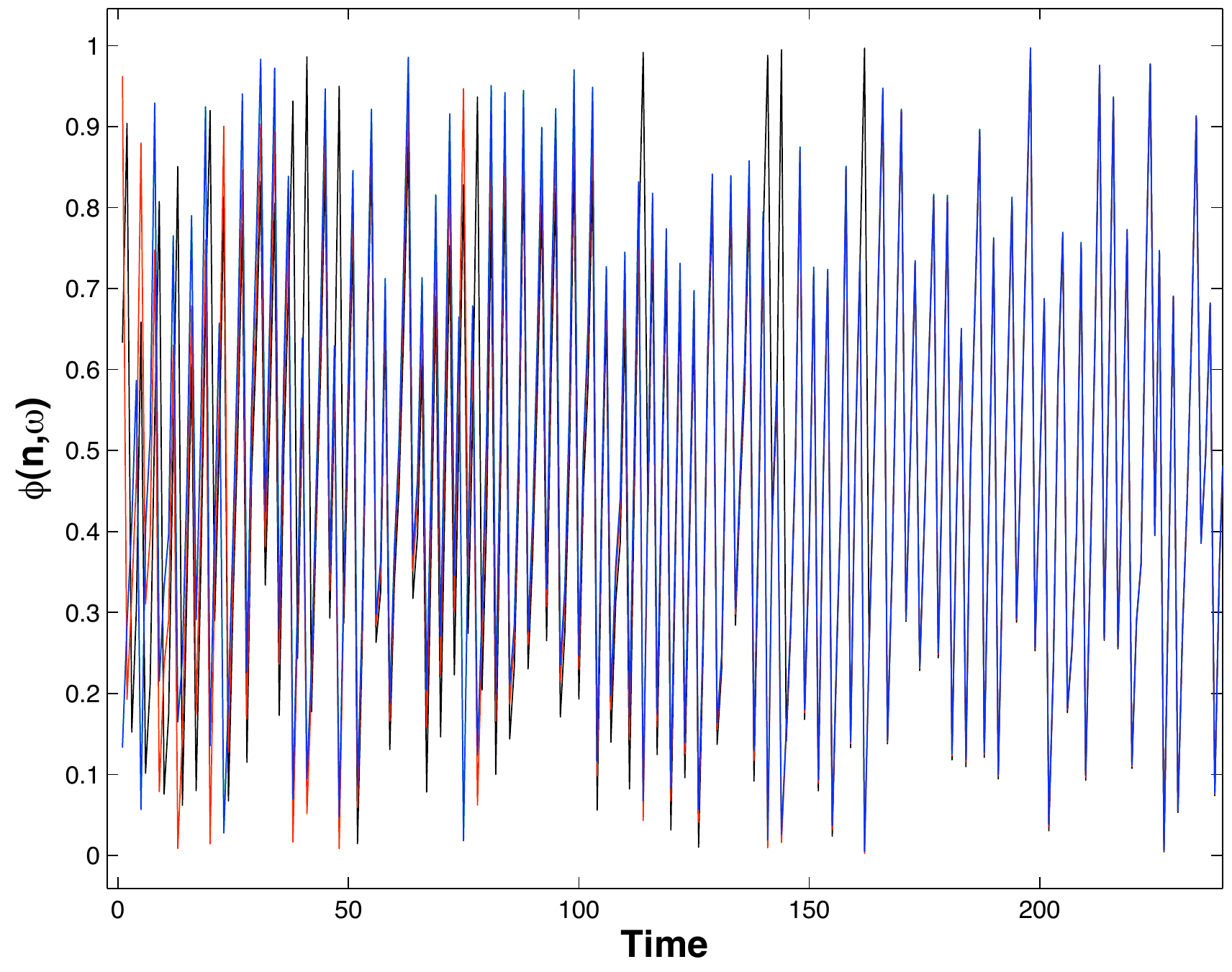}
\caption{Synchronization by additive noise: three distinct
trajectories (in blue, red and black) of
$x_{n+1}=F_{\tau,\epsilon;\omega}(x_n)$, with
$F_{\tau,\epsilon;\omega}$ given by (\ref{Eq_Circle_map_noisy}); the
three trajectories start from three initial points on the circle,
but are driven by the same realization $\omega$ of the noise, and
thus converge to the same random fixed point $a(\omega)$, which is
moving with time. The parameters are $\epsilon = 0.5$, $\tau =
0.283$ and $\sigma = 0.3$, and the corresponding Lyapunov exponent
is $\lambda \simeq -0.0104$.} \label{synchro}
\end{figure}

We illustrate this remarkable property in Fig. \ref{synchro} in the
case of a random fixed point, for given $\epsilon$ and $\tau$. The
Lyapunov exponent for the three distinct trajectories shown in Fig.
\ref{synchro} is strictly negative and the trajectories are
exponentially attracted to the single random fixed point
$a(\omega)$, the realization of the driving system $\theta(\omega)$
being the same for all the trajectories; see Appendices A and B.
Kaijser \cite{Kaijser} provided rigorous results on this type of
synchronization phenomenon, but in a totally different conceptual
setting. Interestingly, as the noise intensity increases, the
Lyapunov exponent becomes more negative, so that the synchronization
occurs even more rapidly, given a fixed realization $\omega$.

This clustering behavior of trajectories with different initial data
is in fact well known for flows on the circle \cite{Kloeden}. In our
example, this phenomenon in phase space is related to a smoothing of
the Devil's staircase in parameter space, the latter cannot be
solely explained by the former. Indeed, we show in Appendix B that
for different irrational numbers and a sufficiently high noise
level, the corresponding stochastic dynamics are stochastically
equivalent, an equivalence that results in the smoothing of certain
steps of the Devil's staircase.

As shown in the lower panel of Fig. \ref{noisytong}, there is also a
direct relationship between the random dynamics and the support of
the PDF on the circle. For a given noise level, this support can
either be the union of a finite number of disjoint intervals (red
and blue curves) or it can fill the whole circle (black curve). The
random attractor is, accordingly, either a random periodic orbit,
with the disjoint intervals being visited in succession, or a random
fixed point; this PDF behavior characterizes the level of the noise
needed to destroy a given tongue.

An exact definition of random fixed point and random periodic orbit
is given in Appendix B, where we provide a rigorous justification of
the numerical results in Figs. \ref{tong}, \ref{noisytong}, and
\ref{synchro}. This theoretical analysis helps clarify the
interaction between noise and nonlinear dynamics in the context of
the GCM classification problem we are interested in.

\section{Concluding remarks}

We recall that Section 2 dealt with the natural, interannual and
interdecadal variability of the ocean's wind-driven circulation. The
oceans' internal variability is an important source of uncertainty
in past-climate reconstructions and future-climate projections
\cite{NRC95, MG0, dijk_ghil, HDbook}. In Section 3 and Appendices A
and B, we dealt more generally with the problem of structural
instability as a possible cause for the stubborn tendency of the
range of uncertainties in climate change projections to increase,
rather than diminish over the last three decades \cite{Charney,
hough91, hough01, SPM}; see again Fig. 1. We summarize here the main
results of the two sections in succession, and outline several open
problems.

The wind-driven double-gyre circulation dominates the near-surface
flow in the oceans'  mid-latitude basins. Particular attention was
paid to the North Atlantic and North Pacific, traversed by the
best-known oceanic jets, namely the Gulf Stream and the Kuroshio
Extension (see Fig. 2). The wind-driven circulation exhibits very
rich internal dynamics and multiscale behavior associated with
turbulent mesoscales (see Fig. 3). Aside from the intrinsic interest
of  this problem in physical oceanography, these major oceanic
currents help regulate the climate of the adjacent continents, while
their low-frequency variability affects past, present and future
global climate.

Thanks in part to the systematic use of dynamical systems theory, a
comprehensive understanding of simple, barotropic, quasi-geostrophic
(QG) models of the double-gyre circulation has been achieved over
the last two decades, and was reviewed in Section 2 here. In
particular, the importance of symmetry-breaking and homoclinic
bifurcations (see Fig. 4) in explaining the observed low-frequency
variability has been validated across a wide hierarchy of models,
including models with much more comprehensive physical formulation,
more realistic geometry, and greater resolution in the horizontal
and vertical \cite{dijk_ghil, HDbook}. This successive-bifurcation
theory also provides a self-consistent explanation for the
climatically important 7-year and 14-year oscillations of the
oceanic circulation and the related atmospheric phenomena in and
around the North-Atlantic basin
\cite{dijk_ghil,HDbook,Mor98,W99,DC&CDV,PGV,Dubar, GK,SGITW2,SGD05}.

The next challenge in physical oceanography is to reconcile the
points of view of dynamical systems theory and statistical mechanics
in describing the interaction between the largest scales of motion
and geostrophic mesoscale turbulence, which is fully captured in
baroclinic QG models. We emphasize that the complexity of these
models of the double-gyre circulation is intermediate between
high-end GCMs and simple ``toy" models; these models offer,
therefore, an ideal laboratory to test our ideas. In particular,
stochastic parametrizations of the rectification process, absent in
barotropic QG models, could be studied using some of the concepts
and tools from RDS theory presented here. Note that the RDS approach
has already been used in the context of stochastic partial
differential equations, in particular for showing the existence of
random attractors, as well as stable, unstable and inertial
manifolds. Thus RDS concepts and tools are not restricted to
finite-dimensional systems \cite{robinson, Duan04, Duan03}.

In Section 3, we have addressed the range-of-uncertainty problem for
IPCC-class GCM simulations (see Fig. 1) by considering them as
stochastically perturbed dynamical systems. This approach is
consonant with recent interest for stochastic parametrizations in
the high-end  modeling-and-simulation community \cite{LN00, LN02,
LN03, TP00, JPS, SZG}. Rigorous mathematical results from the
dynamical systems literature suggest that --- in the absence of
stochastic ingredients --- GCMs as well as simpler models, found on
the lower rungs of the modeling hierarchy \cite{GR00}, are bound to
differ from each other in their results.

This sensitivity follows from the fact that, among deterministic
dynamical systems, those that are hyperbolic are essentially the
only ones that are also structurally stable,  at least in the $C^1$
case \cite{ro,rob, mane, Palis}. Thus, because hyperbolic systems
are not dense in the set of smooth deterministic ones
\cite{Smale66}, we are led to conclude that the topological,
structural-stability approach does not guarantee deterministic-model
robustness, in spite of its many valuable contributions so far.
Related issues for GCM modeling were emphasized recently by Mitchell
\cite{mitchell}, Held \cite{held} and  McWilliams \cite{mcwill}.

We have gone one step further and considered model robustness in the
presence of stochastic terms; such terms could represent either
parametrizations of unresolved processes in GCMs or stochastic
components of natural or anthropogenic forcing, such as volcanic
eruptions or fluctuations in greenhouse gas or aerosol emissions.
Despite the obvious gap between idealized models and high-end
simulations, we have brought to bear random dynamical systems (RDS)
theory \cite{LArnold} on the former.

 In this framework, we have considered a robustness criterion that could
replace structural stability, through the concept of stochastic
conjugacy (see Figs. \ref{RDSbundle} and \ref{RA}). We have shown,
for a stochastically perturbed Arnol'd family of circle maps, that
noise can enhance model robustness. More precisely, this circle map
family exhibits structurally stable, as well as structurally
unstable behavior. When noise is added, the entire family exhibits
{\it stochastic} structural stability, based on the
stochastic-conjugacy concept, even in those regions of parameter
space where deterministic structural instability occurs for
vanishing noise (see Figs. \ref{tong} and \ref{noisytong}).

Clearly the hope that noise can smooth the very highly structured
pattern of distinct  behavior types for climate models, across the
full hierarchy, has to be tempered by a number of caveats. First,
serious questions remain at the fundamental, mathematical level
about the behavior of nonhyperbolic chaotic attractors in the
presence of noise \cite{anischenko_et_al, anischenko_et_al2,
anischenko_et_al3}. Likewise, the case of driving by nonergodic
noise is being actively studied \cite{Arnold_Xu, Arnold_Imkeller,
Li}.

Second, the presence of certain manifestations of a Devil's
staircase has been documented across the full hierarchy of ENSO
models \cite{GR00,EN1,EN2,EN3,EN4,EN6,SG01,EN5}, as well as in
certain observations \cite{GR00, EN5}. Interestingly, both GCMs and
observations only exhibit a few, broad steps of the staircase, such
as $4:1=4$ yr, $4:2=2$ yr, and $4:3\cong 16$ months. Does this
result actually support the idea that nature and its detailed models
always provide sufficient noise to achieve considerable smoothing of
the much finer structure apparent in simpler models? Be that as it
may, we need a much better understanding of how different types of
noise --- additive and multiplicative, white and colored --- act
across even a partial hierarchy of models, say from the simplest
ones, like those studied in Section 3, to the intermediate ones
considered in Section 2.

Third, one needs to connect more closely the nature of a stochastic
parametrization and its effects on the model's behavior in
phase-parameter space. As shown in Appendix B, not all types of
noise are equal with respect to these effects. We are thus left with
a rich, and hopefully fruitful, set of questions, which we expect to
pursue in future work.

\vspace{2ex}

\noindent {\bf Acknowledgements}: It is a pleasure to thank the
organizers and participants of the Conference on the ``{\it Euler
Equations: 250 Years On}" and, more than all, Uriel Frisch, for a
stimulating and altogether pleasant experience. We are grateful to
I. M. Held, J. C. McWilliams, J. D. Neelin and I. Zaliapin for many
useful discussions and their continuing interest in the questions
studied here. A. Sobolevski{\u\i}, E. Tziperman and an anonymous
reviewer have provided constructive criticism and stimulating
remarks that helped improve the presentation. This study was
supported by the U.S. Department of Energy grant DE-FG02-07ER64439
from its Climate Change Prediction Program, and by the European
Commission's No. 12975 (NEST) project ``Extreme Events: Causes and
Consequences (E2-C2)."

\appendix
\section{RDSs and random attractors}
We present here briefly the mathematical concepts and tools of
random dynamical systems, random attractors and stochastic
equivalence. We shall use the concept of pullback attractor
introduced in Section 3.2.1 to define the closely related notion of
a random attractor, but need first to define an RDS. We denote by
$\mathbb{T}$ the set $\mathbb{Z}$, for maps, or $\mathbb{R}$, for
flows. Let $(X,\mathcal{B})$ be a measurable phase space, and
$(\Omega, \mathcal{F}, \mathbb{P}, (\theta(t))_{t\in \mathbb{T}})$
be a {\em metric} dynamical system {\it i.e.} a flow in the
probability space $(\Omega, \mathcal{F}, \mathbb{P})$, such that
$(t,\omega)\mapsto \theta(t) \omega$ is measurable and $\theta(t):
\Omega \rightarrow \Omega$ is measure preserving, $i.e., \theta(t)
\mathbb{P}=\mathbb{P}$.

Let $\varphi: \mathbb{T}\times \Omega \times X \rightarrow X$,
$(t,\omega,x)\mapsto \varphi (t,\omega) x,$ be a mapping with the
two following properties:

\begin{itemize}
\item[] (R$_1$): $\varphi(0,\omega)=\mbox{Id}_X$, and
\item[] (R$_2$) (the cocycle property): For all $s,t \in \mathbb{T}$ and all $\omega \in \Omega$,
$$\varphi(t+s, \omega) =\varphi(t,\theta(s)\omega) \circ \varphi(s,\omega).$$
\end{itemize}
If $\varphi$ is measurable, it is called a {\em measurable} RDS over
$\theta$. If, in addition, $X$ is a topological space (respectively
a Banach space), and $\varphi$ satisfies $(t,\omega)\mapsto
\varphi(t,\omega) x$ continuous (resp. $C^k$, $1\leq k\leq \infty$)
for all $(t,\omega) \in \mathbb{T}\times \Omega$, then $\varphi$ is
called a {\em continuous} (resp. $C^k$) RDS over the flow $\theta$.
If so, then \be (\omega, x)\mapsto \Theta(t)(x,\omega): =
(\theta(t)\omega,\varphi(t,\omega)x),\label{skew} \de is a
(measurable) flow on $\Omega\times X$, and is called the {\em
skew-product} of $\theta$ and $\varphi$. In the sequel, we shall use
the terms ``RDS" or ``cocycle" synonymously.

\begin{figure}[htpb]
\centering
\includegraphics*[width=9cm]{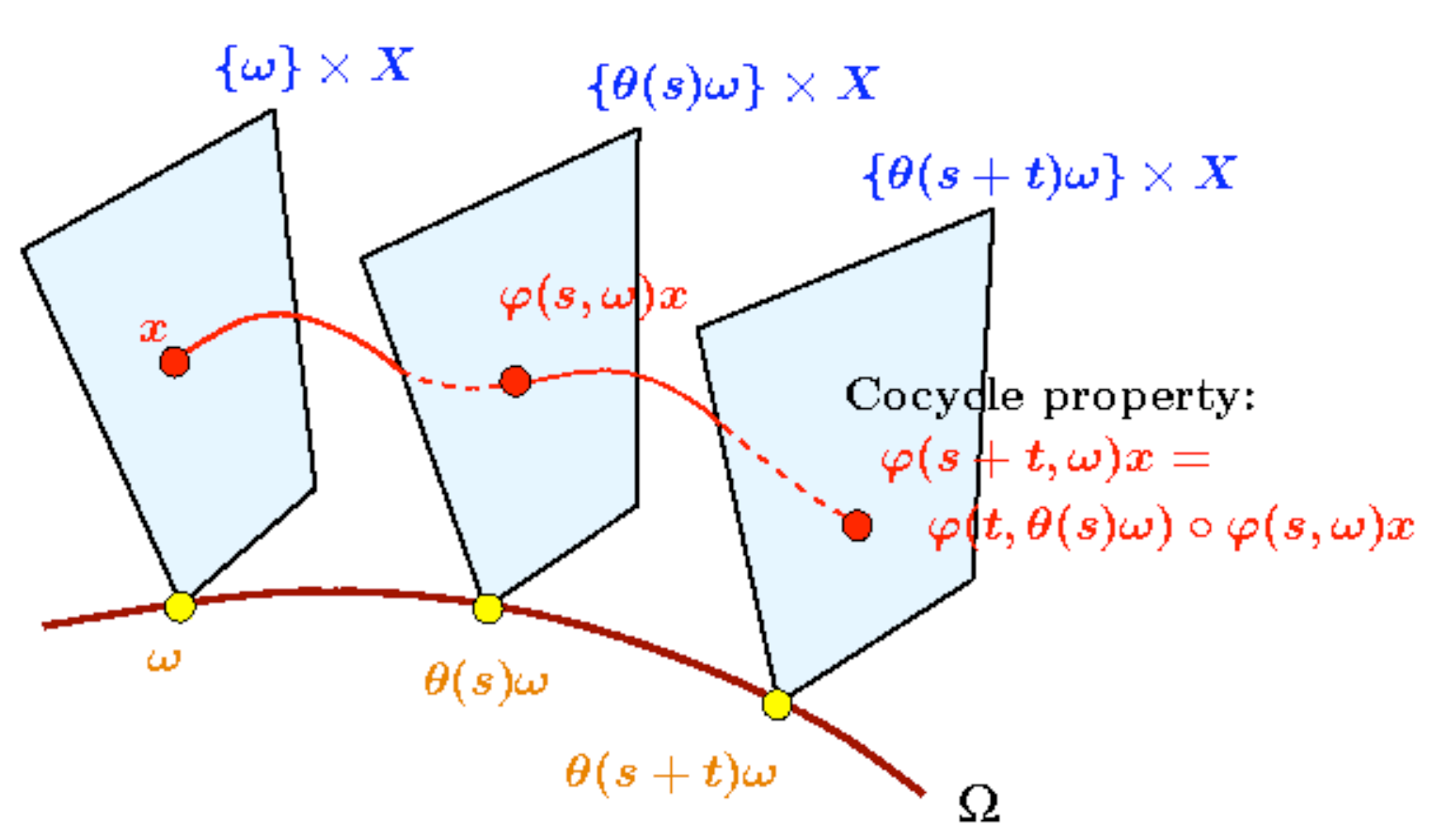}
\caption{Random dynamical systems (RDS) viewed as a flow on the
bundle $X \times \Omega$ = ``dynamical space" $\times$ ``probability
space." For a given state $x$ and realization $\omega$, the RDS
$\varphi$ is such that $\Theta(t)(x,\omega) =
(\theta(t)\omega,\varphi(t,\omega)x)$ is a flow on the bundle.}
\label{RDSbundle}
\end{figure}

The choice of the so-called {\it driving system} $\theta$ is a
crucial step in this set-up; it is mostly dictated by the fact that
the coupling between the stationary driving and the deterministic
dynamics should respect the time invariance of the former, as
illustrated in Fig. \ref{RDSbundle}. The driving system $\theta$
also plays an important role in establishing {\it stochastic
conjugacy} \cite{Cong_book} and hence the kind of classification we
aim at.

\begin{figure}[htpb]
\centering
\includegraphics*[width=9cm]{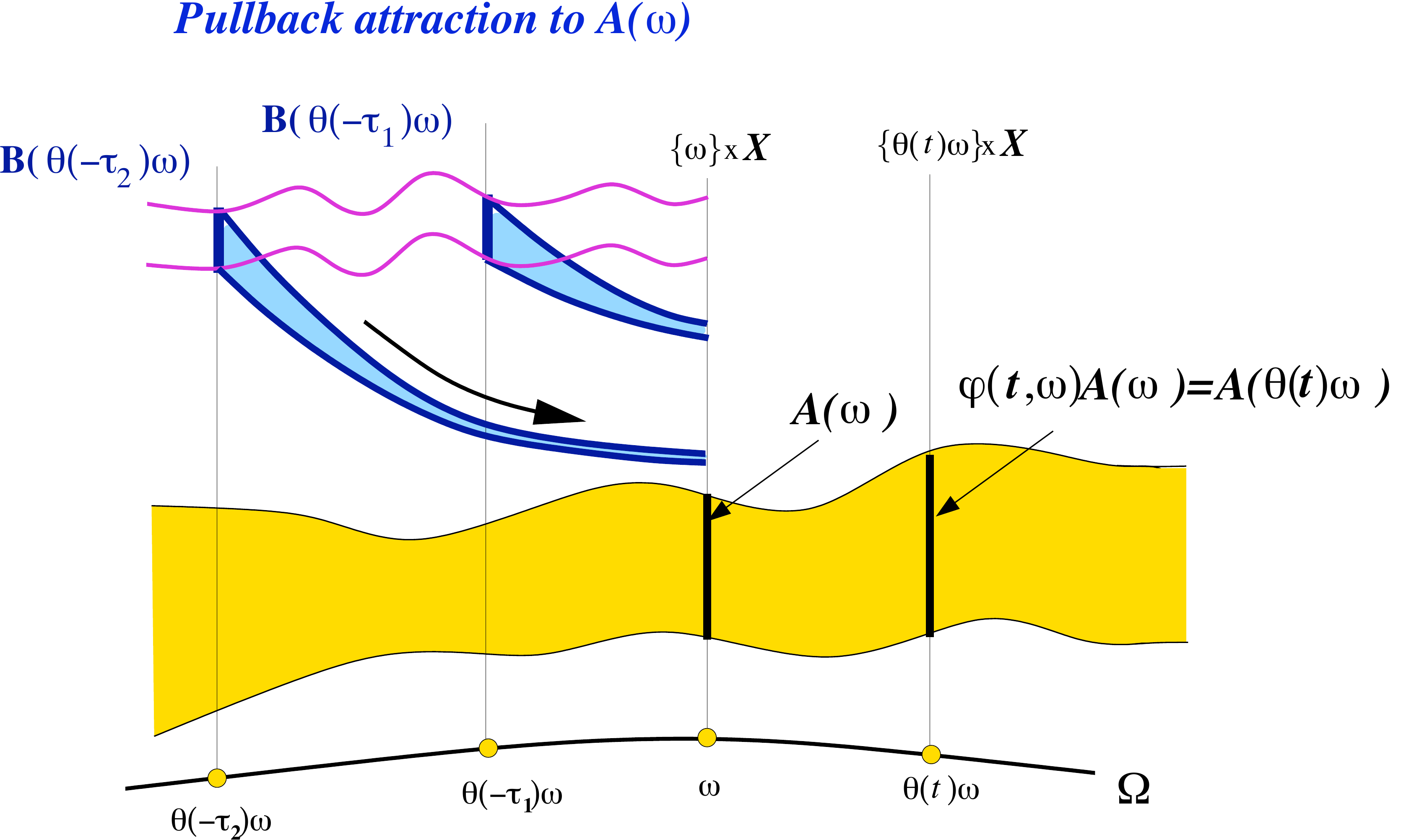}
\caption{Schematic diagram of a random attractor ${\cal A}(\omega)$,
where $\omega \in \Omega$ is a fixed realization of the noise. To be
attracting, for every set $B$ of $X$ in a family $\mathfrak{B}$ of
such sets, one must have $\lim_{t\to +\infty} {\rm
dist}(B(\theta(-t)\omega), {\cal A}(\omega)) = 0$ with
$B(\theta(-t)\omega):=\varphi(t,\theta(-t)\omega) B$; to be
invariant, one must have $\varphi(t,\omega) {\cal A}(\omega) = {\cal
A}(\theta(t)\omega)$. This definition depends strongly on
$\mathfrak{B}$; see \cite{Crauel} for more details.} \label{RA}
\end{figure}

The concept of {\em random attractor} is a natural and
straightforward extension of the definition of pullback attractor
(\ref{attr}), in which Sell's \cite{sell} process is replaced by a
cocycle, {\it cf.} Fig. \ref{RDSbundle}, and the attractor ${\cal
A}$ now depends on the realization $\omega$ of the noise, so that we
have a family of random attractors ${\cal A}(\omega)$, {\it cf.}
Fig. \ref{RA}. Roughly speaking, for a fixed realization of the
noise, one ``rewinds" the noise back to $t\rightarrow -\infty$  and
lets the experiment evolve (forward in time) towards a possibly
attracting set ${\cal A}(\omega)$; the driving system $\theta$
enables one to do this rewinding without changing the statistics,
cf. Figs. \ref{RDSbundle} and \ref{RA}.

Other notions of attractor can be defined in the stochastic context,
in particular based on the original SDE; see \cite{CrauelFlan} or
\cite{Crauel} for a discussion on this topic. The present
definition, though, will serve us well.

Having defined RDSs and random attractors, we now introduce the
notion of {\it stochastic equivalence} or {\it conjugacy}, in order
to rigourously compare two RDSs; it is defined as follows: two
cocycles $\varphi_1(\omega,t)$ and $\varphi_2(\omega,t)$ are
conjugated if and only if there exists a random homeomorphism $h \in
{\rm Homeo}(X)$ and an invariant set such that $h(\omega)(0) = 0$
and \be \label{stochequi} \varphi_1(\omega,t) =
h(\theta(t)\omega)^{-1} \circ \varphi_2(\omega,t) \circ h(\omega).
\de Stochastic equivalence extends classic topological conjugacy to
the bundle space $X \times \Omega$, stating that there exists a
one-to-one, stochastic change of variables that continuously
transforms the phase portrait of one sample system in $X$ into that
of any other such system.

\section{Coarse-graining of the circle map family} \label{Rigorous_section}

We provide here a rigorous justification of the numerical results
obtained in Section \ref{Num_section} on the topological
classification of the family of Arnol'd circle maps in the presence
of noise. Consider the following random family of diffeomorphisms:

\be \label{Eq_Circle_map_noisy} F_{\tau, \epsilon;
\omega}(x):=x+\tau+\sigma \omega -\epsilon \sin(2\pi x) \; \mbox{
mod} \; 1 ,\de for $x\in \mathbb{S}^1$, $\epsilon$ a real parameter
in $(0,1)$, and $\omega$ a random parameter distributed in the
compact interval $I=[-1/2,1/2]$ with fixed distribution $\nu$ and
noise intensity $\sigma$. We denote by $F_{\tau, \epsilon}$ the
corresponding deterministic family of diffeomorphisms when the noise
is switched off, $\sigma=0$.

In the RDS framework, we need to specify the metric dynamical system
modeling the noise. We choose here the interval $\sigma I$ as the
base for the probability space $\Omega$ and define the flow $\theta$
simply as mapping the point $\omega$ into its successor in a
sequence of realizations of the noise. One could also use an
irrational rotation on $\Omega$ for instance; in either case,
ergodicity is ensured.

For the sake of simplicity, we omit for the moment the dependence on
$\tau$ and $\epsilon$. In discrete time, with
$\mathbb{T}=\mathbb{Z}$, we define a map $\phi : \mathbb{T}\times
\Omega \times \mathbb{S}^1 \rightarrow \mathbb{S}^1$, $(n,\omega,
x)\mapsto \phi(n,\omega)x$, such that

\be \phi(n,\omega):= \left \{\begin{array}{l}
F_{\theta^{n-1}\omega}\circ \cdot \cdot \cdot
\circ F_{\omega}, \quad n \geq 1,\\
\mbox{Id}_{\mathbb{S}^1},  \qquad \qquad \qquad n=0,\\
F_{\theta^{n}\omega}^{-1}\circ \cdot \cdot \cdot \circ
F_{\theta^{-1}\omega}^{-1}, \quad n \leq -1.
\end{array}\right.
\de One can prove easily that this $\phi$ satisfies the cocycle
property and is in fact a $C^{\infty}$ RDS on $\mathbb{S}^1$ over
$\theta$.

The pair of mappings $\Theta:= (\theta, \phi)$ is the corresponding
skew-product (\ref{skew}), and it defines a flow on $\Omega \times
\mathbb{S}^1$ by the relation:

\be (\omega,x)\mapsto \Theta(n)(\omega,x):=(\theta^n \omega,
\phi(n,\omega)x).\de A stationary measure $m$ on $\mathbb{S}^1$
under the random diffeomorphism $F_{\tau, \epsilon;\omega}$ yields a
 $\Theta$-invariant measure $\mu:=m\times \nu$, {\it i.e.}
$\Theta_n\mu=\mu$; explicitly,  \be \left.
\begin{array}{l}\int_{\Omega\times \mathbb{S}^1}
f(\omega,x)\mu(d\omega,dx)=\\
\int_{\Omega\times \mathbb{S}^1}
f(\theta_n\omega,\phi(n,\omega)x)\mu(d\omega,dx)\end{array}\right.\de
for all $n \in \mathbb{T}$ and $f\in L^1(\Omega\times \mathbb{S}^1,
\mu)$.

Let us recall the following important proposition \cite{zmarrou}
concerning the stationary measures obtained from the random family
$\{F_{\tau, \epsilon;\omega}\}$.

\begin{theorem}\label{measure_class}

The random circle diffeomorphism $F_{\tau, \epsilon;\omega}$ has a
unique stationary measure $m_{\tau,\epsilon}$. The support of
$m_{\tau,\epsilon}$ consists either of $q$ mutually disjoint
intervals or of the entire circle $\mathbb{S}^1$. The density
function $\phi_{\tau, \epsilon}$ is in $C^{\infty} (\mathbb{S}^1)$
and depends $C^{\infty}$ on $\tau$.  The invariant measure $\mu$ is
ergodic. If the support of $m$ is connected, then it is mixing and
so is $\mu$.
\end{theorem}
Mixing for $m$ means that, for any bounded function
$f:\mathbb{S}^1\rightarrow \mathbb{R}$, and  for an arbitrary
initial point $x_0 \in \mathbb{S}^1$,
$\mathbb{E}(f(\phi(n,\omega)x_0)$ tends to $\int_{\mathbb{S}^1}
f(x)m(dx)$ as $n\rightarrow +\infty$; see \cite{kuksin} for more on
random attractors and mixing.

For deterministic diffeomorphisms of the circle, the rotation number
measures the average rotation per iterate of $F_{\tau, \epsilon}$.
In the presence of noise, one can still define a rotation number for
$F_{\tau, \epsilon; \omega}$, namely \be \rho_{\tau,
\epsilon;\omega}(x)=\lim_{k \to \infty} \frac{\tilde{F}^k_{\tau,
\epsilon; \omega} (x)-x}{k}, \de where $\tilde{F}$ denotes the lift
of a map $F$, acting on $\mathbb{S}^1$ modulo $1$, to a map acting
on $\mathbb{R}$. For fixed $\tau$ and $\epsilon$, we can then show
that $\rho_{\tau, \epsilon;\omega}$ exists for $\nu$-almost all
$\omega$ and is a constant; this constant $\rho_{\tau;\omega}$ is
independent of $x$ and $\omega$ \cite{zmarrou}. Furthermore, $\tau
\rightarrow \rho_{\tau, \epsilon}$ is $C^{\infty}$ for each
$\epsilon$, which is not true in the deterministic case with $\sigma
= 0$; see again \cite{zmarrou}.

Theorem \ref{measure_class} has a natural geometric counterpart in
terms of random attractors, as confirmed through our numerical
study; see again Fig. \ref{noisytong}. More precisely, we introduce
also the following definitions of {\it random fixed point} and {\it
random periodic orbit}; these definitions differ somewhat from those
given in \cite{zmarrou}.

\begin{definition}\label{random_fixed_point}
A random fixed point is a measurable map $a : \Omega \rightarrow
\mathbb{S}^1$ for which \be
\phi(1,\omega)a(\omega)=a(\theta(\omega)), \de for $\nu$-almost all
$\omega \in \Omega,$ {\it i.e.} such that $\Omega \times a(\Omega)$
is an invariant set for the flow  given by the skew-product
$\Theta$. A random periodic orbit of period $q$ is likewise an
invariant set with cardinality $q$ in fibers $\mathbb{S}^1\times
\{\omega\}$ for $\nu$-almost all $\omega$.
\end{definition}

With these definitions, the following results of \cite{zmarrou}
still hold.

\begin{theorem} \label{dyn_behav}
For a random diffeomorphism $F_{\tau,\epsilon; \omega}$ of the
circle $\mathbb{S}^1$, with a stationary measure $ m$ supported on a
union $E$ of $q$ disjoint intervals, the corresponding skew-product
$\Theta$ restricted to $E$ has precisely one attracting random
periodic orbit and one repelling random periodic orbit.
\end{theorem}
Attraction in the preceding theorem means that $\lim_{n\rightarrow
\infty} |F_{\omega}^n(x)-F_{\omega}^n(a(\omega))|=0$, for a set of
initial data $(x,\omega)\in \mathbb{S}^1\times \Omega$ with positive
$\lambda \times \nu$-measure, in the case of a random attracting
fixed point; here $\lambda$ is Lebesgue measure on $\mathbb{S}^1$
and the extension to a random periodic orbit is obvious.

Using these two theorems and rigorous results on random point
attractors \cite{Crauel}, we can show that (i) if the support of the
stationary measure is the whole circle (black curve in Fig.
\ref{noisytong}), then there exists one random fixed point which is
pullback attracting; and (ii) if the support consists of $q$
disjoint intervals, then the random attractor is a random periodic
orbit of period $q$ (red and blue curves).

Having explained how the connectedness of the PDF support at
different noise levels is related to the nature of the random
attractor, we now turn to an explanation of the ``disappearance" of
the smaller steps in the Devil's staircase, as the noise level
increases. To do so, we consider the Lyapunov spectrum of an RDS,
which still relies on the Oseledets \cite{Oseledets} {\em
multiplicative ergodic theorem} (MET).

To state an MET for RDS on manifolds, we differentiate
$\phi(n,\omega)$ at $x\in \mathbb{S}^1$, and obtain the linear map
\be T\phi (n,\omega,x): T_x M \rightarrow T_{\phi(n,\omega)x} M, \de
where $T \phi$ is a continuous linear cocycle on the tangent bundle
$TM$ of the manifold $M$ over the skew-product flow $\Theta$. If the
flow $\phi$ possesses an ergodic invariant measure $\mu$ such that
the required integrability condition for applying the MET is
verified with respect to $\mu$, then the MET holds for $\phi$ over
$M$ \cite{arnold_L_trend}.

Because of Theorem \ref{measure_class} here, we can apply the MET to
our problem and conclude that a unique Lyapunov exponent exists for
the {\em linearization} of each diffeomorphism belonging to our
family of random diffeomorphisms, and that this exponent is
independent of the realization of the noise. We show next how to use
the Lyapunov spectrum in studying the stochastic equivalence classes
of a given RDS family, along with its  driving system $\theta$. This
last aspect of the classification problem is outlined for linear
hyperbolic cocycles.

N.D. Cong \cite{Cong} has shown that, even in the linear context,
the main difference with respect to the deterministic case is that
the classification depends strongly on the properties of $\theta$,
which is directly linked to the system noise and its modeling. For
instance, if $\theta$ is an irrational rotation on $\mathbb{S}^1$,
one can construct infinitely many classes of hyperbolic cocycles
that are not pairwise topologically equivalent, by playing
essentially with the orientations of the cocyles, {\it i.e.}
reversing between clockwise and anticlockwise rotation on
$\mathbb{S}^1.$ As we shall see, such difficulties can be avoided in
the case of noisy Arnol'd tongues, especially for additive noise.
Related issues still form an active research area in RDS theory; see
\cite{arnold_L_trend} for a brief survey.

A key ingredient for the linear classification is the notion of {\em
coboundary}, which we recall herewith.

\begin{definition}
A measurable set $K \subset \Omega$ is called a {\em coboundary} if
there exists a set $H\in \mathcal{F}$ such that $K=H\triangle \theta
H$, where $H\triangle \theta H$ denotes the symmetric difference of
$H$ and $\theta H$.
\end{definition}

Let $A$ and $B$ be two linear random maps on $\mathbb{R}^d$, and
denote by deg $A(\omega)$ and deg $B(\omega)$ the degrees of the
maps $A(\omega)$ and  $B(\omega)$ with respect to a chosen random
orientation. These degrees are just the sign of the determinant of
the corresponding random matrices, and equal $-1$ or $1$; see
\cite{Cong_book, Cong} for details. Consider the two linear
hyperbolic cocycles $\Phi_A$ and $\Phi_B$, associated with the maps
$A$ and $B$, and the following subset of $\Omega$:

\be C_{AB}=\{\omega \in \Omega|\mbox{deg}\; A(\omega)\cdot
\mbox{deg} \;B(\omega)=-1\}; \de $C_{AB}$ is just the set of all
$\omega \in \Omega$ for which the degrees of the two linear maps
$A(\omega)$ and $B(\omega)$ differ.

The main theorem for the classification of our diffeomorphisms of
$\mathbb{S}^1$ follows \cite{Cong}.

\begin{theorem} \label{Cong_classification}
Two one-dimensional linear hyperbolic cocycles $\Phi_A$ and $\Phi_B$
are conjugate if and only if the following conditions hold:
\begin{itemize}
\item[] (i) $\mbox{sign} \; \lambda_A =\mbox{sign} \; \lambda_B$, and
\item[] (ii) the associated set $ C_{AB}$ is a coboundary.
\end{itemize}
\end{theorem}
Here $\lambda_{A}$ and $\lambda_{B}$ indicate the Lyapunov exponents
of $\Phi_A$ and $\Phi_B$, respectively.

Before applying this result, let us explain heuristically how a
Devil's staircase step that corresponds to a rational rotation
number can be ``destroyed" by a sufficiently intense noise. Consider
the period-1 locked state in the deterministic setting. At the
beginning of this step, a pair of fixed points is created, one
stable and the other unstable. As the bifurcation parameter is
increased, these two points move away from each other, until they
are $\pi$ radians apart. Increasing the parameter further causes the
fixed points to continue moving along, until  they finally meet
again and are annihilated in a {\em saddle-node bifurcation},  thus
signaling the end of the locking interval.

When noise is added, we have to distinguish between a ``strongly
locked" regime, where the stable and unstable fixed points are
nearly $\pi$ radians apart, and a ``weakly locked"  regime, where
these two fixed points are close to each other. In both regimes, the
relaxation time in the vicinity of the stable point represents an
important time scale of the problem. In the strongly locked regime,
this is the only time scale of interest. In the weakly locked
regime, though, the process of escaping across the unstable fixed
point is nonnegligible and the associated escape time becomes the
second time scale of interest. From these heuristic considerations
it follows that the distinction between strong and weak locking
depends on the strength of the external noise.

If we consider period-$T$ locked states, with $T\geq 2$, the same
kind of reasoning can be applied to the stable and unstable
$T$-cycle. We conclude therefore, for a fixed  $\epsilon > 0$, that
the narrower Devil's staircase steps are the least robust, while the
wider ones are the most robust.

The fact that a locked case becomes unlocked  when noise is growing
implies in particular that the rotation number
$\rho_{\tau,\epsilon}$ becomes {\em irrational}  for a sufficiently
high noise level. According to Theorem 2.1 of \cite{Kaijser}, the
Lyapunov exponent is strictly {\em negative} in this case almost
surely. Moreover, by reinterpreting other results of Kaijser
\cite{Kaijser} in our RDS framework, we can show that the random
attractor is in fact a random fixed point; this, in turn, allows us
to conclude that the corresponding linearized cocycle at this random
fixed point is hyperbolic. Next, by using the Hartman-Grobman
theorem for RDSs \cite{wanner,ruffino, HGT}, we can conjugate the
nonlinear cocycle with its linearization; in fact, Theorem 3.1 of
\cite{ruffino} says that this conjugacy is global.

Consider now two linearized cocycles $\Phi_A$ and $\Phi_B$, at one
and the same or at two distinct random fixed points of the family of
random diffeomorphisms, for the same noise intensity, and denote by
$A(\omega)$ and $B(\omega)$ the random linear parts of the cocycles
$\Phi_A$ and $\Phi_B$ respectively; it follows from our model of
noisy circle maps that $C_{AB}$ is empty. Indeed, the noise being
additive, the random orientation is preserved for different
parameter values. But $\theta$ is assumed to be ergodic, and so we
have that $\Omega\triangle \theta \Omega $ is empty and, therefore,
$C_{AB}$ is a coboundary. Therewith, Theorem
\ref{Cong_classification} can be applied to obtain the desired
result for the problem considered here: with an appropriate amount
of noise, two deterministic diffeomorphisms that are not
topologically equivalent can fall into the same topological
stochastic class! The numerical results  of Section 3.2 are entirely
in agreement with this assertion.

Note that the set $C_{AB}$ could differ from a coboundary, if the
noise occurred additively in the phase of the nonlinear term, for
instance. Here we see the importance of noise modeling in obtaining
the conjectural view of Fig. \ref{amib} for a family of dynamical
systems in general.

It follows, in particular, that the exact nature of the stochastic
parametrizations in a family of GCMs does matter. It's not enough to
follow the trend by devising and implementing such parametrizations:
one should test that a given parametrization, once found to be
suitable in other respects, does improve the proximity, in an
appropriate sense, between climate simulations within the family of
GCMs for which it has been been developed.


\vspace{1ex}

\end{document}